\documentclass[12pt,openany,letterpaper,twoside]{book}

\usepackage{amsmath,amssymb,amscd}
\usepackage{amsthm}
\usepackage{stmaryrd}
\usepackage{fancyhdr}
\usepackage{times}
\usepackage{mathrsfs} 
\usepackage{titlesec}
\usepackage{remreset}
\usepackage[pdftex,dvips]{geometry}

\input diagxy

\newtheoremstyle{plain}
     {\topsep}
     {\topsep}
     {\itshape}
     {}
     {\bfseries}
     {}
     { }
     {\thmnumber{#2.}\hspace{0.5ex}\thmname{#1.}\thmnote{ \rm (#3)}}

\newtheoremstyle{definition}
     {\topsep}
     {\topsep}
     {\normalfont}
     {}
     {\bfseries}
     {}
     { }
     {\thmnumber{#2.}\hspace{0.5ex}\thmname{#1.}\thmnote{ \rm (#3)}}

\newtheoremstyle{fact}
     {\topsep}
     {\topsep}
     {\slshape}
     {}
     {\bfseries}
     {}
     { }
     {\thmnumber{#2.}\hspace{0.5ex}\thmname{#1.}\thmnote{ \rm (#3)}}

\theoremstyle{plain}
\newtheorem{theorem}[subsection]{Theorem}
\newtheorem*{theorem*}{Theorem} 
\newtheorem{lemma}[subsection]{Lemma}
\newtheorem{proposition}[subsection]{Proposition}
\newtheorem{corollary}[subsection]{Corollary}

\theoremstyle{definition}

\newtheorem{remarks}[subsection]{Remarks}
\newtheorem{remark}[subsection]{Remark}
\newtheorem*{remark*}{Remark}

\newtheorem*{question*}{Question}
\newtheorem*{examples*}{Examples}  

\newtheorem*{example*}{Example}

\newtheorem*{convention*}{Convention}

\theoremstyle{fact}

\newtheorem{ftheorem}[subsection]{Theorem}
\newtheorem{flemma}[subsection]{Lemma}
\newtheorem{fproposition}[subsection]{Proposition}
\newtheorem{fcorollary}[subsection]{Corollary}

\def\proofont{\fontseries{bx}\selectfont}
\def\proofname{Proof. }

\newcommand{\pcite}[2]{{\cite[#1]{#2}}}

\newcommand{\rest}{\mbox{\parbox[t]{0.1cm}{$|$\\[-10pt] $|$}}}

\newcommand{\Note}[1]{}

\makeatletter
\renewenvironment{proof}[1][\proofname]{\par
  \normalfont
  \topsep6\p@\@plus6\p@ \trivlist
  \item[\hskip\labelsep\noindent\proofont #1]\ignorespaces
}{%
  \qed\endtrivlist
}
\makeatother

\titlelabel{\thetitle.\ }
\titleformat*{\subsection}{\normalsize\bfseries}
\titlespacing{\subsection}{0pt}{\topsep}{0.5ex}
\titleformat{\subsection}[runin]{\normalfont\bfseries}{%
\thesubsection.}{0.5ex}{}[.]

\makeatletter
  \renewcommand{\thechapter}{\Roman{chapter}}
  \renewcommand{\thesection}{\thechapter.\Alph{section}}
  \renewcommand{\thesubsection}{\thechapter.\arabic{subsection}}
  \@removefromreset{subsection}{section}
  \@addtoreset{subsection}{chapter}
\makeatother

\author{G\'abor Luk\'acs
\thanks{I gratefully acknowledge the generous financial support received
from the Killam Trusts and Dalhousie University that enabled me to do 
this research.}}
   
\title{Notes on duality theories of abelian groups
}

\begin{document}

\makeatletter
\let\mytitle\@title
\chead{\small\itshape G. Luk\'acs / \mytitle }
\fancyhead[RO,LE]{\small \thepage}
\makeatother

\maketitle

\def\thanks#1{} 

\thispagestyle{empty}


\newpage

\pagenumbering{roman}
\pagestyle{plain}
\setcounter{page}{2}

\mbox{ }

\newpage

\renewcommand{\thepage}{(\roman{page})}
\tableofcontents

\newpage
\thispagestyle{empty}

\newpage

\mbox{ }

\newpage

\pagenumbering{arabic}
\setcounter{page}{1}

\pagestyle{fancy}

\chapter*{Introduction}
\addcontentsline{toc}{chapter}{\protect {}Introduction}

\thispagestyle{empty}

In December 2005, I had the good fortune of spending a week as the guest
of Reinhard Winkler at Technische Univert\"at Wien, and had the honour of
being asked to deliver a mini-course on the topic of duality theory of
abelian groups. During this visit, I came to realize that the only
elementary text on this topic that I could refer a graduate student to,
namely the one by Pontryagin \cite{Pontr}, was published 20 years ago, and
its well-organized content is half a century old. This is not to say that
there are no excellent books that contain a lot of results on abelian
topological groups (such as \cite{DikProSto}, \cite{Rudin}, \cite{Morris},
or \cite{Aussenhofer}), but my feeling is that they depart from the aim of
presenting basic notions of duality theory in a self-contained and
elementary manner. Therefore, encouraged by my host and his research group
in Vienna, I decided to write up the notes that I prepared for the mini
course (Chapters~1 and~2), and to keep developing it as a longer term
project. In order to provide a smoother presentation, some elementary
results of were collected in the Appendix.

I use {\itshape italics font} for results that appear to be new and not 
part of the "common knowledge," while {\slshape slanted font} is used for 
all other statements. Let me know if I am wrong about any of them.

\section*{Past and future work}

This notebook is under development, so all suggestions, comments, 
and questions are warmly welcome.

\subsection*{Chapter 1}
This is definitely the only more-or-less ready piece, but I am not happy
with my treatment of locally compact abelian groups. (Similarly to
Roeder's approach (cf.~\cite{Roeder3}), the only result borrowed from
functional analysis is the Peter-Weyl theorem. I am still missing a nice
proof of the structure theorem of compactly generated LCA groups, and a
complete proof of the Pontryagin duality.)

\subsection*{Chapter 2}
I have only a sketch. After a categorical introduction, which explains 
why cartesian closed categories are so interesting, I hope to cover 
k-groups (of Noble) and convergence groups in this chapter.

\subsection*{Appendix} 
It evolves as the chapters develop. It is a collection of results many 
readers might be familiar with, so I saw no point in including them into 
the chapters.

\subsection*{Future topics} Precompact (abelian) groups [the problem is 
that it requires the reader to be familiar with notions such as $C$-embedded 
subsets, etc]; examples of pathological or otherwise interesting groups 
[although \cite{DikProSto} is probably the best source for this]; localic 
groups [it would be a nice example of group objects in a category].

\subsection*{Topics left out} Nuclear groups, topological vector spaces 
[these are both important, but I am not sure if I can present them in an 
elementary way].

\section*{Acknowledgements}

I am deeply indebted to Reinhardt Winkler and his research group for their 
hospitality during my visit at Technische Univert\"at Wien, who inspired 
and encouraged me to organize my thoughts and notes on the topic.

I wish to thank Salvador Hern\'andez, Keith Johnson, and Jeff Egger for
the valuable discussions and helpful suggestions that were of great
assistance in writing this paper.

\newpage

\thispagestyle{empty}

\chapter{The Pontryagin dual}

\thispagestyle{empty}

\section{The evaluation homomorphism}

\nocite{Roeder3}
\nocite{Roeder2}

\subsection{Pontryagin dual}

For Hausdorff topological groups $G$ and $H$, we denote by 
$\mathscr{H}(G,H)$ the space of continuous homomorphisms 
$\varphi\colon G\rightarrow H$ equipped with the compact-open topology.
(For a brief review of the compact-open topology, see Appendix~\ref{app:A}.)
Since the property of being a homomorphism is equationally defined,
$\mathscr{H}(G,H)$ is a closed subspace of $\mathscr{C}(G,H)$. Indeed, 
\begin{align}
\mathscr{H}(G,H) = \bigcap\limits_{g_1,g_2\in G} 
\{f \in \mathscr{C}(G,H) \mid f(g_1 g_2)=f(g_1)f(g_2)\},
\end{align}
and each of the sets in the intersection is closed.
Topological groups have natural uniform structures, and the compact-open 
topology can also be realized as the topology of uniform convergence on 
compacta: $\varphi_\alpha \longrightarrow \varphi$ in $\mathscr{H}(G,H)$ 
if $\varphi_\alpha\rest_K$ converges uniformly to $\varphi\rest_K$ 
for every compact subset $K \subseteq G$.

For $A\in \mathsf{Ab(Top)}$, a {\em character} of $A$ is a homomorphism 
$\chi\colon A \rightarrow \mathbb{T}$, where 
$\mathbb{T}:=\mathbb{R}/\mathbb{Z}$. If $\chi$ is continuous, then it must 
factor through the maximal Hausdorff quotient $A/N_A$ of $A$, because 
$\mathbb{T}$ is Hausdorff (cf.~Proposition~\ref{prop:max:T2}). The {\em 
Pontryagin dual} $\hat A$ of $A$ is the group 
$\mathscr{H}(A/N_A,\mathbb{T})$ of all continuous characters of $A$ 
equipped with the compact-open topology. Since compact subsets of $A/N_A$ 
are precisely the images of compact subsets in $A$, one may view $\hat A$ 
as the group $\mathscr{H}(A,\mathbb{T})$ both algebraically and 
topologically (cf.~Corollary~\ref{cor:max:T2}).

We put $A_d$ for the group $A$ equipped with the discrete topology. The 
space $\mathscr{C}(A_d,\mathbb{T})$ coincides with $\mathbb{T}^{A_d}$, and 
thus compact. Therefore, its closed subspace $\widehat{A_d}$ is also 
compact; it carries the topology of {\em pointwise convergence} on $A$, in 
other words, $\chi_\alpha \rightarrow \chi$ if and only if $\chi_\alpha(a) 
\rightarrow\chi(a)$ for every $a \in A$. Since the compact-open topology 
is finer than the pointwise convergence one, $\hat A \rightarrow 
\widehat{A_d}$ is continuous.

\subsection{Polar} The sets $\Lambda_n=[-\frac 1 {4n}, \frac 1 {4n}]$ form 
a base at $0$ for $\mathbb{T}$, and have the property that if
$k x \in \Lambda_1$ for $k=1,\ldots,n$, then $x \in \Lambda_n$.
(Here, $[- \frac  1 {4n}, \frac 1 {4n}]$ are identified with their
images in $\mathbb{T}$.)
For $S\subseteq A$ and $\Phi \subseteq \hat A$, their 
{\em polar} sets are defined as
\begin{align}
S^\vartriangleright & = \{ \chi \in \hat A \mid \chi(S)\subseteq 
\Lambda_1\}, \\
\Phi^\vartriangleleft & = \{ a \in A \mid \forall \chi \in \Phi,
\chi(a) \in \Lambda_1 \} = 
\bigcap\limits_{\chi \in \Phi} \chi^{-1}(\Lambda_1).
\end{align}
Since $S^\vartriangleright$ is closed in the topology of pointwise 
convergence, in particular, it is closed in the compact-open one. 
On the other hand, 
$\Phi^\vartriangleleft$ is closed in $A$ because each $\chi \in\Phi$ is 
continuous. It is immediately seen that 
$(\vartriangleleft,\vartriangleright)$ is a Galois-connection between 
subsets of $A$ and $\hat A$. Thus, $S \subseteq 
S^{\vartriangleright\vartriangleleft}$ and 
$\Phi\subseteq \Phi^{\vartriangleleft\vartriangleright}$.

\begin{lemma} \label{lemma:pol:hom}
Let $\varphi\colon A \rightarrow B$ be a morphism in $\mathsf{Ab(Top)}$.

\begin{list}{{\rm (\alph{enumi})}}
{\usecounter{enumi}\setlength{\labelwidth}{25pt}\setlength{\topsep}{-8pt} 
\setlength{\itemsep}{-5pt} \setlength{\leftmargin}{25pt}}
 
\item For $S \subseteq A$,
$\hat\varphi^{-1}(S^\vartriangleright)=\varphi(S)^\vartriangleright$.

\item
For $\Sigma \subseteq \hat B$,
$\varphi^{-1}(\Sigma^\vartriangleleft)=\hat\varphi(\Sigma)^\vartriangleleft$.

\item
For $S^\prime\subseteq B$, 
$\hat\varphi(S^{\prime\vartriangleright}) \subseteq 
\varphi^{-1}(S^\prime)^\vartriangleright$.

\item
For $\Sigma^\prime \subseteq \hat A$, 
$\varphi(\Sigma^{\prime\vartriangleleft}) \subseteq
\hat\varphi^{-1}(\Sigma^\prime)^\vartriangleleft$.

\item
For  $S \subseteq A$, $\varphi(S^{\vartriangleright\vartriangleleft})
\subseteq \varphi(S)^{\vartriangleright\vartriangleleft}$.

\item
For $S^\prime\subseteq B$, 
$\varphi^{-1}(S^\prime)^{\vartriangleright\vartriangleleft} \subseteq
\varphi^{-1}(S^{\prime\vartriangleright\vartriangleleft})$.

\end{list}
\end{lemma}

\begin{proof} 
Since $\varphi$ is continuous, it induces a continuous homomorphism
$\hat\varphi\colon \hat B \rightarrow \hat A$ defined by
$\hat\varphi(\chi)= \chi \circ \varphi$.
\begin{align}
\chi \in \varphi(S)^\vartriangleright & \Longleftrightarrow 
\chi(\varphi(S))\subseteq \Lambda_1  \Longleftrightarrow 
\chi\circ\varphi \in S^\vartriangleright  \Longleftrightarrow
\chi \in \hat\varphi^{-1}(S^\vartriangleright),
\\
\intertext{and so (a) follows.}
g \in \varphi^{-1}(\Sigma^\vartriangleleft) & \Longleftrightarrow
\varphi(g) \in \Sigma^\vartriangleleft \Longleftrightarrow 
\forall \chi \in \Sigma, \chi(\varphi(g))\in\Lambda_1  \\
& \Longleftrightarrow
\forall \psi \in \hat\varphi(\Sigma), \psi(g) \in \Lambda_1
\Longleftrightarrow g \in \hat\varphi(\Sigma)^\vartriangleleft,
\end{align}
which shows (b).

(c) Since $\varphi(\varphi^{-1}(S^\prime)) \subseteq S^\prime$, one has
$S^{\prime\vartriangleright} \subseteq 
\varphi(\varphi^{-1}(S^\prime))^\vartriangleright$. By (a) applied to 
$S=\varphi^{-1}(S^\prime)$, 
$\varphi(\varphi^{-1}(S^\prime))^\vartriangleright=
\hat\varphi^{-1}(\varphi^{-1}(S^\prime)^\vartriangleright)$. Thus, 
$ S^{\prime\vartriangleright} \subseteq 
\hat\varphi^{-1}(\varphi^{-1}(S^\prime)^\vartriangleright)$.

(d) Since $\hat\varphi(\hat\varphi^{-1}(\Sigma^\prime)) \subseteq 
\Sigma^\prime$, one has 
$\Sigma^{\prime\vartriangleleft}\subseteq 
\hat\varphi(\hat\varphi^{-1}(\Sigma^\prime))^\vartriangleleft$.
By (b) applied to $\Sigma=\hat\varphi^{-1}(\Sigma^\prime)$, 
$\hat\varphi(\hat\varphi^{-1}(\Sigma^\prime))^\vartriangleleft = 
\varphi^{-1}(\hat\varphi^{-1}(\Sigma^\prime)^\vartriangleleft)$.
Thus,  $\Sigma^{\prime\vartriangleleft}\subseteq 
\varphi^{-1}(\hat\varphi^{-1}(\Sigma^\prime)^\vartriangleleft)$.

(e) follows from (a) and (d), and (f) follows from (b) and (c).
\end{proof}

\begin{flemma} \label{lemma:sub:ann}
Let $A \in \mathsf{Ab(Top)}$, $H \leq A$ its subgroup, and let 
$\pi_H\colon A \rightarrow A/H$ be the canonical projection.
Then:

\begin{list}{{\rm (\alph{enumi})}}   
{\usecounter{enumi}\setlength{\labelwidth}{25pt}\setlength{\topsep}{-8pt}
\setlength{\itemsep}{-5pt} \setlength{\leftmargin}{25pt}}

\item
$H^\vartriangleright$ is a subgroup;

\item
$\hat\pi_H\colon \widehat{A/H} \rightarrow \hat A$ is injective, and its 
image is $H^\vartriangleright$;

\item 
{\rm (\pcite{5}{Brug2})}
if $S^\prime \subseteq A/H$ and $0 \in S^\prime$, then
$\pi_H^{-1}(S^\prime)^\vartriangleright=\hat\pi_H(S^{\prime\vartriangleright})$.

\end{list}
\end{flemma}

\begin{proof}
(a) Since $\Lambda_1$ contains only the zero subgroup,
$\chi(H)\subseteq \Lambda_1$ holds if and only if
$H \subseteq \ker \chi$. Thus, $H^\vartriangleright$ is the annihilator 
subgroup of $H$ in $\hat A$.

(b) Since $\pi_H$ is onto, $\hat\pi_H$ is injective. 
Clearly, $\operatorname{Im} \hat \pi_H \subseteq H^\vartriangleright$, 
so in order to show the converse, let $\chi \in H^\vartriangleright$.
Then $H \subseteq \ker \chi$, and thus $\chi$ induces a continuous 
character $\bar\chi\colon A/H \rightarrow \mathbb{T}$ such that
$\hat\pi_H(\bar \chi)=\bar \chi \circ  \pi_H=\chi$.

(c) One has  $H \subseteq \pi_H^{-1}(S^\prime)$ because $0 \in S^\prime$, 
and thus $\pi_H^{-1}(S^\prime)^\vartriangleright \subseteq 
H^\vartriangleright=\operatorname{Im} \hat\pi_H$, by (b).
Lemma~\ref{lemma:pol:hom}(c) applied to $S=\pi_H^{-1}(S^\prime)$ yields
$\hat\pi_H^{-1}(\pi_H^{-1}(S^\prime)^\vartriangleright)=
\pi_H(\pi_H^{-1}(S^\prime))^\vartriangleright=S^{\prime\vartriangleright}$.
This completes the proof, because $\hat \pi_H$ is injective.
\end{proof}

For $X \in \mathsf{Top}$, we denote by $\mathcal{K}(X)$ the collection of 
compact subsets of $X$. For $G \in \mathsf{Grp(Top)}$, one 
puts $\mathcal{N}(G)$ for the collection of neighborhoods of $e$ in $G$.

\begin{samepage}

\begin{fproposition} \label{prop:pol:basic}
Let $A \in \mathsf{Ab(Top)}$. 

\begin{list}{{\rm (\alph{enumi})}}
{\usecounter{enumi}\setlength{\labelwidth}{25pt}\setlength{\topsep}{-8pt} 
\setlength{\itemsep}{-5pt} \setlength{\leftmargin}{25pt}}

\item
$\chi \in \hom(A,\mathbb{T})$ is continuous if and 
only if there is $U \in \mathcal{N}(A)$ such that 
$\chi(U)\subseteq \Lambda_1$.

\item
The collection $\{K^\vartriangleright\mid K\subseteq 
A,K\in\mathcal{K}(A)\}$ is a base at $0$ to $\hat A$.

\item
$\Sigma \subseteq \hat A$ is equicontinuous if and 
only if there is $U \in \mathcal{N}(A)$ such that 
$\Sigma\subseteq U^\vartriangleright$.

\end{list}
\end{fproposition}

\end{samepage}

\begin{proof}
Observe that in order to establish continuity-like properties of
homomorphisms of topological groups, it suffices to check them 
at a single point, and $0$ is 
a convenient choice for that purpose. 

(a) Necessity is obvious. So, in order to show sufficiency, 
let $n \in \mathbb{N}$ and let $U \in \mathcal{N}(A)$ be such that
$\chi(U)\subseteq \Lambda_1$. Because of the continuity of 
addition in $A$, there is $V \in \mathcal{N}(A)$ such that
$\underbrace{V + \cdots + V}_\text{n times} \subseteq U$, and in 
particular,  $k\chi(V)\subseteq \Lambda_1$ for every $k=1,\ldots,n$.
Therefore, $\chi(V)\subseteq \Lambda_n$, as desired.

(b) Let $C \in \mathcal{K}(A)$, and set $K=C\cup 2C\cup\cdots\cup nC$.
Then $K \in \mathcal{K}(A)$, and $\chi(C)\subseteq \Lambda_n$
for every $\chi \in K^\vartriangleright$,  because 
$k\chi(a) \in \Lambda_1$ for every $a \in C$ and  $k=1,\ldots,n$.
Thus, we obtained $K^\vartriangleright \subseteq
\{\chi \in \hat A \mid \chi(C)\subseteq \Lambda_n\}$, as desired.

(c) Because of the homogeneous structure of topological groups,
equicontinuity of a family $\Sigma$ of characters can be checked at a 
single point. Equicontinuity at $0$ means that for every neighborhood 
$\Lambda_n$ of $0$ there is $V \in \mathcal{N}(A)$ such that
$\chi(V)\subseteq \Lambda_n$ for every $\chi \in \Sigma$. By the argument 
presented in (a), this condition is satisfied for all $n$ if and only if 
it is satisfied for $n=1$. Since this condition for $n=1$ is precisely
$\Sigma \subseteq U^\vartriangleright$, the proof is complete.
\end{proof}

\subsection{Evaluation homomorphism}
Each $a \in A$ gives rise to a continuous character $\hat a$ of
$\hom(A,\mathbb{T})$ given by evaluation: $\hat a(\psi)=\psi(a)$. In 
particular, $\hat a$ is a continuous character of $\hat A$ (whose 
topology is finer than that of $\hom(A,\mathbb{T})$). Thus,
$\hat a \in \hat{\hat A}$, and the  evaluation map
\begin{align}
\alpha_A\colon A & \longrightarrow \hat{\hat A} \\
a & \longmapsto \hat a
\end{align}
is a homomorphism of groups. 
In the sequel, we present necessary and sufficient conditions for 
$\alpha_A$ to be  continuous and an embedding (cf.
Propositions~\ref{prop:ev:cont} and~\ref{prop:ev:emb}), and it will 
become clear that in general, $\alpha_A$ need not be continuous.

A collection $\mathcal{C}$ of compact sets 
of a topological space $X$ is a {\em cobase} if for every compact subset 
$K \subseteq X$ there is $C \in \mathcal{C}$ such that $K \subseteq C$.

\begin{samepage}

\begin{fproposition}[\pcite{2.3}{Noble}] \label{prop:ev:cont}
Let $A \in \mathsf{Ab(Top)}$.

\begin{list}{{\rm (\alph{enumi})}}
{\usecounter{enumi}\setlength{\labelwidth}{25pt}\setlength{\topsep}{-8pt}
\setlength{\itemsep}{-5pt} \setlength{\leftmargin}{25pt}}

\item
For every $U \in \mathcal{N}(A)$, $U^\vartriangleright$ is compact in 
$\hat A$.

\item
The following statements are equivalent:

\begin{list}{{\rm (\roman{enumii})}}
{\usecounter{enumii}\setlength{\labelwidth}{25pt}\setlength{\topsep}{-6pt}
\setlength{\itemsep}{-2pt} \setlength{\leftmargin}{25pt}}

\item
$\alpha_A$ is continuous;

\item
every compact subset of $\hat A$ is equicontinuous;

\item
$\{U^\vartriangleright \mid U \in \mathcal{N}(A)\}$ is a cobase for
$\hat A$;

\item
$\{U^{\vartriangleright\vartriangleright} \mid U \in \mathcal{N}(A)\}$ is 
a base for $\hat{\hat A}$ at $0$.

\end{list}

\item
$\ker \alpha_A=\bigcap\limits_{\chi \in \hat A} \ker \chi$.
\end{list}

\end{fproposition}

\end{samepage}

\begin{proof}
(a) As noted earlier,  $U^\vartriangleright$ is closed in $\hat A$ 
(it is closed even in $\hom(A,\mathbb{T})$).
By Proposition~\ref{prop:pol:basic}(c), it is also equicontinuous, and
therefore it is compact by a standard Arzela-Ascoli type argument.

(b) The equivalence of (ii) and (iii) is an immediate consequence of 
Proposition~\ref{prop:pol:basic}(c). 

(i) $\Rightarrow$ (ii): Suppose that $\alpha_A$ is continuous, and let
$\Sigma \subseteq \hat A$ be compact. Then $\Sigma^\vartriangleright$ is 
a neighborhood of $0$ in $\hat{\hat A}$. Thus, by continuity of 
$\alpha_A$, there is $U \in 
\mathcal{N}(A)$ such that $\alpha_A(U) \subseteq\Sigma^\vartriangleright$,
which is equivalent to $\Sigma \subseteq U^\vartriangleright$. Therefore,
by Proposition~\ref{prop:pol:basic}(c), $\Sigma$ is equicontinuous.

(iii) $\Rightarrow$ (iv):
If $V \in \mathcal{N}(\hat{\hat A})$, then by 
Proposition~\ref{prop:pol:basic}, there is a compact subset
$\Phi\subseteq \hat A$ such that $\Phi^\vartriangleright \subseteq V$. 
Since $\{U^\vartriangleright \mid U \in \mathcal{N}(A)\}$ is a cobase for 
$\hat A$, there is $U \in \mathcal{N}(A)$ such that 
$\Phi \subseteq U^\vartriangleright$. Therefore, 
$U^{\vartriangleright\vartriangleright} \subseteq \Phi^\vartriangleright 
\subseteq V$, as desired.

(iv) $\Rightarrow$ (i): If $W\in \mathcal{N}(\hat{\hat A})$, then there is 
$U \in \mathcal{N}(A)$ such that 
$U^{\vartriangleright\vartriangleright} \subseteq W$, and therefore
$\alpha_A(U) \subseteq U^{\vartriangleright\vartriangleright}\subseteq W$.
Hence, $\alpha_A$ is continuous.

(c) For $a \in A$, one has
\begin{align}
a \in \ker \alpha_A \Longleftrightarrow \alpha_A(a)=0 \Longleftrightarrow 
\forall \chi \in \hat A, (\alpha_A(a))(\chi)=0 \Longleftrightarrow
\forall \chi \in \hat A, \chi(a)=0,
\end{align}
which completes the proof.
\end{proof}

A combination of Propositions~\ref{prop:pol:basic}(b) 
and~\ref{prop:ev:cont}(b) yields:

\begin{fcorollary} \label{cor:ev:conthat}
Let $A \in \mathsf{Ab(Top)}$. The following statements are equivalent:

\begin{list}{{\rm (\roman{enumi})}}
{\usecounter{enumi}\setlength{\labelwidth}{25pt}\setlength{\topsep}{-8pt}
\setlength{\itemsep}{-4pt} \setlength{\leftmargin}{25pt}}

\item
$\alpha_{\hat A}$ is continuous;

\item
$\{K^{\vartriangleright\vartriangleright} \mid K \in\mathcal{K}(A)\}$
is a cobase for $\hat{\hat A}$. \qed
\end{list}
\end{fcorollary}

A map $f\colon X \rightarrow Y$ between Hausdorff spaces is said to be 
{\em $k$-continuous} if the restriction $f\rest_K$ is continuous for every 
compact subset $K$ of $X$. Although $\alpha_A$ need not be continuous 
(cf.~Proposition~\ref{prop:ev:cont}), its restriction to any compact subset of 
$A$ is continuous.

\begin{ftheorem} \label{thm:ev:kcont}
Let $A \in \mathsf{Ab(Haus)}$. The evaluation homomorphism $\alpha_A$ is 
$k$-continuous.
\end{ftheorem}


\begin{proof}
By Corollary~\ref{cor:CO:kcont}, $e\colon A\rightarrow 
\mathscr{C}(\mathscr{C}(A,\mathbb{T}),\mathbb{T})$ is $k$-continuous. 
Since the image $e(A)$ is contained in 
$\mathscr{H}(\mathscr{C}(A,\mathbb{T}),\mathbb{T})$, and
the map $\mathscr{H}(\mathscr{C}(A,\mathbb{T}),\mathbb{T}) \rightarrow
\mathscr{H}(\mathscr{H}(A,\mathbb{T}),\mathbb{T})=\hat{\hat A}$ (given by 
restriction to $\mathscr{H}(A,\mathbb{T})$) is continuous, 
$\alpha_A$ is the composition of a $k$-continuous map with a continuous 
one. 
\end{proof}

\subsection{Local quasi-convexity}
Let $A \in \mathsf{Ab(Top)}$.
Following Banaszczyk, $S\subseteq A$ is {\em quasi-convex} if 
$S=S^{\vartriangleright\vartriangleleft}$. The group $A$ is
{\em locally quasi-convex} (or briefly, {\em LQC}) if it admits a base of 
quasi-convex neighborhoods at $0$, that is, if
$\{U^{\vartriangleright\vartriangleleft} \mid U \in \mathcal{N}(A)\}$ 
is a base at $0$ (cf. \cite{Banasz}). 

\begin{fproposition}[\pcite{1}{Brug2}] \label{prop:qc:hat}
The group $\hat A$ is locally quasi-convex for every
$A \in \mathsf{Ab(Top)}$.
\end{fproposition}

\begin{proof}
By Proposition~\ref{prop:pol:basic}, the collection $\{K^\vartriangleright 
\mid K \in \mathcal{K}(A)\}$ is
a base at $0$ at $\hat A$. Each member of the collection is quasi-convex, 
because $K^{\vartriangleright\vartriangleright\vartriangleleft}=
K^\vartriangleright$. 
\end{proof}

A variant of the next proposition appeared in \cite[6.10]{Aussenhofer} and 
\cite[4.3]{Butz2}, and seems to be a well-known result.

\begin{fproposition} \label{prop:ev:emb}
Let $A \in \mathsf{Ab(Top)}$.

\begin{list}{{\rm (\alph{enumi})}}
{\usecounter{enumi}\setlength{\labelwidth}{25pt}\setlength{\topsep}{-8pt} 
\setlength{\itemsep}{-5pt} \setlength{\leftmargin}{25pt}}

\item
If $\alpha_A$ is an embedding, then $A$ is locally quasi-convex.

\item
If $A$ is locally quasi-convex, then $\alpha_A$ is open onto 
its image. 

\item
If $A$ is locally quasi-convex and Hausdorff, then $\alpha_A$ is 
injective.

\end{list}
\end{fproposition}

\begin{proof}
First, note that for $S\subseteq A$, one has
$\alpha_A(S^{\vartriangleright\vartriangleleft})=
S^{\vartriangleright\vartriangleright} \cap \alpha_A(A)$. If
$U \in \mathcal{N}(A)$, then $U^\vartriangleright$ is compact in 
$\hat A$ (cf. Proposition~\ref{prop:ev:cont}(a)), and so
$U^{\vartriangleright\vartriangleright}$ is open in $\hat{\hat A}$ by 
Proposition~\ref{prop:pol:basic}(b). Therefore,
$\alpha_A(U^{\vartriangleright\vartriangleleft})$ is open in $\alpha_A(A)$.

(a) Since $\alpha_A$ is an embedding, $\alpha_A(V)$ is open in 
$\alpha_A(A)$  for every $V \in \mathcal{N}(A)$. Thus, by 
Proposition~\ref{prop:ev:cont}(b)(iv), 
there is $U \in \mathcal{N}(A)$ such that 
$\alpha_A(U^{\vartriangleright\vartriangleleft})=
U^{\vartriangleright\vartriangleright}\cap\alpha_A (A) 
\subseteq \alpha_A(V)$, which implies 
$U^{\vartriangleright\vartriangleleft} \subseteq V$, because 
$\alpha_A$ is injective.

(b) Let $V \in \mathcal{N}(A)$, and using LQC pick $U \in \mathcal{N}(A)$ 
such that $U^{\vartriangleright\vartriangleleft} \subseteq V$. Then one 
has $\alpha_A(U^{\vartriangleright\vartriangleleft}) \subseteq \alpha_A(V)$, 
and $\alpha_A(U^{\vartriangleright\vartriangleleft})$ is open in 
$\alpha_A(A)$.

(c) Let $a \in A$ be a non-zero element. Since $A$ is Hausdorff and LQC,  
there exists $U \in \mathcal{N}(A)$ such 
that $a \not \in U^{\vartriangleright\vartriangleleft}$. Thus,
there is $\chi \in  U^{\vartriangleright}$ such that
$(\alpha_A(a))(\chi)=\chi(a) \not\in \Lambda_1$. Therefore, 
$\alpha_A(a) \neq 0$. Hence, $\alpha_A$ is injective. 
\end{proof}

\begin{remarks} \label{rem:MAP}
(1) In Proposition~\ref{prop:ev:emb}(b), the map $\alpha_A$ need not be an 
embedding, because it is not necessarily continuous or injective.

(2) Following Neumann and Wigner, whenever $\alpha_A$ is 
injective, $A$ is said to be {\em 
maximally  almost-periodic} (or briefly, {\em MAP}; cf. \cite{NeuWig}).

(3) For every $A \in \mathsf{Ab(Top)}$, $\alpha_{\hat A}$ is injective. 
Indeed, if $\alpha_{\hat A}(\chi)=0$, then for every $\xi\in \hat{\hat A}$,
$\xi (\chi)=(\alpha_{\hat A}(\chi))(\xi)=0$. Thus, for $a \in A$ and
$\xi=\alpha_A(a)$, one obtains
$\chi(a)=(\alpha_A(a))(\chi)=0$. Therefore, $\chi=0$, as desired.
\end{remarks}

\begin{corollary} \label{cor:ev:emb}
Let $A \in \mathsf{Ab(Haus)}$ be such that $\alpha_A$ is continuous. Then
$\alpha_A$ is an embedding if and only if $A$ is locally quasi-convex.
\qed
\end{corollary}

We denote by $\mathsf{LQC}$ the full subcategory of $\mathsf{Ab(Top)}$ 
formed by the locally quasi-convex groups, and present a functorial method 
of ``turning" every group into an LQC group.

\begin{theorem} \label{thm:qc:epiref}
$\mathsf{LQC}$ is an epireflective subcategory of $\mathsf{Ab(Top)}$, and 
the reflection $A^\vartriangle$ is given by equipping the group 
of $A$ with a new group topology whose base at $0$ is 
$\{U^{\vartriangleright\vartriangleleft}\mid U \in \mathcal{N}(A)\}$. 
\end{theorem}

In less categorical language, Theorem~\ref{thm:qc:epiref} states that 
for every $A \in \mathsf{Ab(Top)}$:

\begin{list}{{\rm (\arabic{enumi})}}
{\usecounter{enumi}\setlength{\labelwidth}{25pt}\setlength{\topsep}{-8pt} 
\setlength{\itemsep}{-5pt} \setlength{\leftmargin}{25pt}}

\item
$A^\vartriangle$ is a topological group;

\item
there is a continuous surjective homomorphism
$\nu_A\colon A \rightarrow A^\vartriangle$ that is natural in $A$;

\item
$\varphi^\vartriangle \colon A^\vartriangle
\rightarrow B^\vartriangle$ is continuous whenever $\varphi\colon A
\rightarrow B$ is a continuous homomorphism;

\item
Every continuous homomorphism $\gamma\colon A \rightarrow C$
into an LQC group $C$ factors uniquely through $\nu_A$, that is,
there is a unique $\gamma^\prime\colon A^\vartriangle \rightarrow C$ such 
that $\gamma=\gamma^\prime \circ \nu_A$.
\end{list}

\begin{proof} (1) Each set of the form 
$U^{\vartriangleright\vartriangleleft}$ is symmetric (i.e., 
$-U^{\vartriangleright\vartriangleleft}=U^{\vartriangleright\vartriangleleft}$),
and it is easily seen that $(U\cap V)^{\vartriangleright\vartriangleleft} 
\subseteq U^{\vartriangleright\vartriangleleft} \cap 
V^{\vartriangleright\vartriangleleft}$.
Let $U \in \mathcal{N}(A)$, 
and let $V \in \mathcal{N}(A)$ be such  that $V+V \subseteq U$. Then
$V^{\vartriangleright\vartriangleleft}+V^{\vartriangleright\vartriangleleft} 
\subseteq 
(V+V)^{\vartriangleright\vartriangleleft}
\subseteq U^{\vartriangleright\vartriangleleft}$. Therefore, 
the proposed candidate for $\mathcal{N}(A^\vartriangle)$ defines a group 
topology on the underlying group of $A$, as desired.

(2) The topology of $A^\vartriangle$ is coarser than that of $A$, so the 
identity homomorphism $A \rightarrow A^\vartriangle$ is continuous (and 
it is obviously natural in $A$).


(3) Let $U^{\vartriangleright\vartriangleleft}\in\mathcal{N}(B^\vartriangle)$.
By Lemma~\ref{lemma:pol:hom}(f), 
$\varphi^{-1}(U^{\vartriangleright\vartriangleleft}) \supseteq
\varphi^{-1}(U)^{\vartriangleright\vartriangleleft}$. Since 
$\varphi$ is continuous, $\varphi^{-1}(U) \in \mathcal{N}(A)$, and 
therefore $\varphi^{-1}(U)^{\vartriangleright\vartriangleleft} \in
\mathcal{N}(A^\vartriangle)$. Hence, $\varphi^\vartriangle$ is continuous.

(4) Since $A$ and $A^\vartriangle$ have the same underlying set, 
uniqueness of $\gamma^\prime$ is clear. If $C \in \mathsf{LQC}$, then
$C=C^\vartriangle$, and therefore $\gamma^\prime=\gamma^\vartriangle$, 
by (3).
\end{proof}

\begin{corollary}
The limit of a family of locally quasi-convex groups formed in 
$\mathsf{Ab(Top)}$ is locally quasi-convex, and coincides with the limit 
formed in  $\mathsf{LQC}$. In particular,  $\mathsf{LQC}$ is closed under 
the formation of arbitrary products and subgroups.
\end{corollary}

\begin{proof}
The first statement is a well-known category theoretical property of 
reflective subcategories (cf.~\cite[IV.3, V.5]{MacLane}), and it implies 
the second one, as products are limits. For the third statement, 
let $A$ be an LQC group and $H \leq A$ be its subgroup. The reflection 
$(A/H)^\vartriangle$ has the same underlying group as $A/H$, and 
therefore
\begin{equation}
H=\operatorname{Eq}(A \two^{\pi_H^\vartriangle}_{0} (A/H)^\vartriangle),
\end{equation}
where $\pi_H\colon A \rightarrow A/H$ is the canonical projection. Since 
equalizers are limits, this completes the proof.
\end{proof}

\begin{corollary}
For every $A\in \mathsf{Ab(Top)}$, the underlying groups of 
$\hat A$ and $\widehat{A^\vartriangle}$ coincide.
\end{corollary}

\begin{proof}
By Proposition~\ref{prop:qc:hat} the group $\hat{\mathbb{Z}}=\mathbb{T}$ 
is LQC. Thus, by Theorem~\ref{thm:qc:epiref},
every continuous character $\chi\colon A \rightarrow \mathbb{T}$ gives 
rise to a continuous $\chi^\vartriangle\colon A^\vartriangle \rightarrow
\mathbb{T}$. Therefore, $\hat A \subseteq \widehat{A^\vartriangle}$ as 
sets. The reverse inclusion is obvious, because  $A^\vartriangle$  carries a 
coarser topology than $A$.
\end{proof}

\begin{proposition}
Let $A \in \mathsf{Ab(Top)}$ such that $\alpha_A$ is continuous. Then:

\begin{list}{{\rm (\alph{enumi})}}
{\usecounter{enumi}\setlength{\labelwidth}{25pt}\setlength{\topsep}{-8pt} 
\setlength{\itemsep}{-5pt} \setlength{\leftmargin}{25pt}}

\item
$\hat A$ and $\widehat{A^\vartriangle}$ have the same compact subsets.

\item
$A^\vartriangle$ coincides with the group $A$ equipped with the initial
topology induced by $\alpha_A$;

\item
$N_{A^\vartriangle}=\ker \alpha_{A^\vartriangle} = \ker \alpha_A$;

\item
$\alpha_{A^\vartriangle}$ is continuous and open onto its image;

\item
$\hat{\hat \nu}_A\colon \hat{\hat A} \rightarrow 
\widehat{\widehat{A^\vartriangle}}$ is an embedding;

\end{list}
\end{proposition}

\begin{proof}
Since $\hat \nu_A\colon \widehat{A^\vartriangle}\rightarrow \hat A$ is 
continuous, every compact subset $\Phi$ of $\widehat{A^\vartriangle}$ is 
also compact in $\hat A$, and thus, by 
Proposition~\ref{prop:ev:cont}(b), there exists $U \in 
\mathcal{N}(A)$ such that $\Phi \subseteq U^\vartriangleright$.

(a) Suppose that $\Sigma \subseteq \hat A$ is compact. Then, by 
Proposition~\ref{prop:ev:cont}(b), 
$\Sigma \subseteq U^\vartriangleright$ for some $U \in
\mathcal{N}(A)$. By Proposition~\ref{prop:ev:cont}(a), 
$U^\vartriangleright=
(U^{\vartriangleright\vartriangleleft})^\vartriangleright$ 
is compact in 
$\widehat{A^\vartriangle}$. Since $\Sigma=\hat\nu_A^{-1}(\Sigma)$ is a 
closed subset of $U^\vartriangleright$, $\Sigma$ is compact
in $\widehat{A^\vartriangle}$.

(b) By Proposition~\ref{prop:ev:cont}(b), 
$\{U^{\vartriangleright\vartriangleright} \mid U \in \mathcal{N}(A)\}$ is 
a base for $\hat{\hat A}$ at $0$, because $\alpha_A$ is continuous.
The statements follows from 
$\alpha_A^{-1}(U^{\vartriangleright\vartriangleright})=
U^{\vartriangleright\vartriangleleft}$.

(c) The first equality follows from (b) and Proposition~\ref{prop:max:T2},
because $\hat{\hat A}$ is Hausdorff. For the second equality, observe that 
by Proposition~\ref{prop:ev:cont}(c), 
$\ker\alpha_A=\bigcap\limits_{\chi\in \hat A} \ker \chi$, and 
by (a), $A$ and $A^\vartriangle$ have the same continuous characters.

(d) By (a), polars of quasi-convex neighborhoods form a cobase to 
$\widehat{A^\vartriangle}$, and thus $\alpha_{A^\vartriangle}$ is 
continuous, by Proposition~\ref{prop:ev:cont}(b). Therefore, 
the statement follows by Proposition~\ref{prop:ev:emb}(b).

(e) By Proposition~\ref{prop:ev:cont}(b), the collections
$\{U^{\vartriangleright\vartriangleright} \mid U \in \mathcal{N}(A)\}$
and $\{V^{\vartriangleright\vartriangleright} \mid V \in 
\mathcal{N}(A^\vartriangle)\}$ are bases at $0$ to 
$\hat{\hat A}$ and $\widehat{\widehat{A^\vartriangle}}$, respectively.
Since every $V \in \mathcal{N}(A^\vartriangle)$ has the form of 
$U^{\vartriangleright\vartriangleleft}$, this completes the proof.
\end{proof}

\begin{theorem} \label{thm:dense:sbgrp}
Let $A \in \mathsf{Ab(Top)}$, and $D \leq A$ be a dense subgroup. Then:

\begin{list}{{\rm (\alph{enumi})}}
{\usecounter{enumi}\setlength{\labelwidth}{25pt}\setlength{\topsep}{-8pt} 
\setlength{\itemsep}{-5pt} \setlength{\leftmargin}{25pt}}

\item
the underlying groups of $\hat D$ and $\hat A$ coincide;

\item
if $\alpha_D$ is continuous, $\hat D$ and $\hat A$ have the same compact 
subsets, and $\alpha_A$ is continuous;

\item
if $D$ is locally quasi-convex, then so is $A$.
\end{list}
\end{theorem}


\begin{proof}
(a) The inclusion $\iota_D\colon D \rightarrow A$ induces a continuous 
homomorphism $\hat \iota_D\colon \hat A \rightarrow \hat D$, which is 
injective, because $D$ is dense in $A$. In order to show that
$\hat\iota_D$ is surjective, let $\chi \in \hat D$. 
For each $a \in A$, there is a net $\{x_\alpha\}\subseteq D$
such that $x_\alpha \longrightarrow a$. 
Thus, $x_\alpha - x_\beta \longrightarrow 0$,
and so $\chi(x_\alpha) - \chi(x_\beta) \longrightarrow 0$. In other words, 
$\{\chi(x_\alpha)\}$ is a Cauchy net in the complete group $\mathbb{T}$, 
and therefore $\lim \chi(x_\alpha)$ exists.
Set $\bar\chi(a)=\lim \chi(x_\alpha)$; it is 
easily seen that $\bar \chi$ is well-defined and continuous. Hence,
$\hat \iota_D(\bar \chi)=\chi$, as desired.

(b) It suffices to show that $\hat \iota_D^{-1}(\Phi)$ is compact in 
$\hat A$ for every compact subset $\Phi \subseteq \hat D$, because 
$\hat\iota_D$ is continuous. By Proposition~\ref{prop:ev:cont}(b), 
since $\alpha_D$ is continuous, $\Phi$ is equicontinuous on the 
dense subgroup $D$, which means that there is $U \in \mathcal{N}(A)$ such 
$\Phi \subseteq (U\cap D)^\vartriangleright$ (cf. 
Proposition~\ref{prop:pol:basic}(b)). It follows from the density of $D$ 
that $\operatorname{Int} U \subseteq \overline{D\cap U}$, and therefore
\begin{equation}
\Phi\subseteq (U\cap D)^\vartriangleright = 
(\overline{U\cap D})^\vartriangleright \subseteq 
(\operatorname{Int} U)^\vartriangleright.
\end{equation}
Hence, $\hat \iota_D^{-1}(\Phi)$ is a closed subset of the compact subset 
$(\operatorname{Int} U)^\vartriangleright$ of $\hat A$ (cf. 
Proposition~\ref{prop:pol:basic}(a)). Since $\operatorname{Int U} \in 
\mathcal{N}(A)$, this also shows that every compact subset of $\hat A$ is 
contained in  $U^\vartriangleright$ for some $U \in \mathcal{N}(A)$, which 
means, by Proposition~\ref{prop:ev:cont}(b), that $\alpha_A$ is 
continuous.

(c) Let $V \in\mathcal{N}(A)$, and pick $V_1 \in \mathcal{N}(A)$ such that
$V_1 - V_1 \subseteq V$, so $\bar V_1 \subseteq \operatorname{Int} V$ (cf. 
Proposition~\ref{prop:sep:topgrp}(a)).
There exists $W_1 \in \mathcal{N}(D)$  such that 
$W_1^{\vartriangleright\vartriangleleft}\cap D\subseteq V_1\cap D$,
because $V_1 \cap D \in \mathcal{N}(D)$ and $D$ is LQC. 
(For $\Sigma\subseteq \hat D$, $\Sigma^\vartriangleleft$ in $D$
is $\hat\iota_D^{-1}(\Sigma)^\vartriangleleft\cap D$, because $\hat D=\hat 
A$ as sets.) There is $W_2 \in \mathcal{N}(A)$ such that $W_1=W_2 \cap D$, 
and since $D$ is dense, one has $W_1^\vartriangleright=W_2^\vartriangleright$.
Thus, $W_1^{\vartriangleright\vartriangleleft} = 
W_2^{\vartriangleright\vartriangleleft}$, and therefore
\begin{align}
\operatorname{Int} (W_2^{\vartriangleright\vartriangleleft}) \subseteq
\overline{W_2^{\vartriangleright\vartriangleleft}\cap D} = 
\overline{W_1^{\vartriangleright\vartriangleleft}\cap D} \subseteq
\overline{V_1\cap D} = 
\bar V_1 \subseteq V.
\end{align}
For $U \in \mathcal{N}(A)$ such that $U + U \subseteq W_2$, 
one has $U^{\vartriangleright\vartriangleleft}+
U^{\vartriangleright\vartriangleleft} \subseteq 
W_2^{\vartriangleright\vartriangleleft}$, and so
$U^{\vartriangleright\vartriangleleft} \subseteq 
\operatorname{Int} (W_2^{\vartriangleright\vartriangleleft})$ (cf.
Proposition~\ref{prop:sep:topgrp}(a)). Hence,
we found $U \in \mathcal{N}(A)$ such that
$U^{\vartriangleright\vartriangleleft} \subseteq V$, as desired.
\end{proof}

\begin{corollary}
The completion of every locally quasi-convex group is locally 
quasi-convex. \qed
\end{corollary}

\subsection{Precompactness} A subset $S\subseteq G$ of $G \in 
\mathsf{Grp(Top)}$ is said to be {\em precompact} if for every $U \in 
\mathcal{N}(A)$ there is a finite subset $F \subseteq G$ such that $S 
\subseteq FU$. Similarly to compactness, if $S,S^\prime \subseteq G$ are 
precompact and $S_1\subseteq S$, then so are $S_1$, $\bar S$, $S +S^\prime$, 
and $\varphi(S)$ for every continuous homomorphism $\varphi\colon G 
\rightarrow H$. Our interest in this property arises from the following 
theorem on uniform spaces (cf. \cite{Isb} and \cite[8.3.16]{Engel6}):

\begin{ftheorem} \label{thm:precomp:cmpct}
A uniform space $(X,\mathscr{U})$ is compact if and only if it is complete 
and precompact.
\end{ftheorem}


\begin{proposition} \label{prop:lqc:precom}
Let $A \in \mathsf{Ab(Haus)}$ be a locally quasi-convex group. Then:

\begin{list}{{\rm (\alph{enumi})}}
{\usecounter{enumi}\setlength{\labelwidth}{25pt}\setlength{\topsep}{-8pt} 
\setlength{\itemsep}{-5pt} \setlength{\leftmargin}{25pt}}

\item
$\alpha_A^{-1}(E)$ is precompact in $A$
for every equicontinuous subset $E\subseteq \hat{\hat A}$;


\item
$K^{\vartriangleright\vartriangleleft}$ is precompact in $A$
for every compact subset $K\subseteq A$;

\item
if $\alpha_{\hat A}$ is continuous, then $\alpha_A^{-1}(F)$ is precompact 
for every compact $F \subseteq \hat{\hat A}$.

\end{list}
\end{proposition}

\begin{proof} 
First, we note that by Proposition~\ref{prop:ev:emb}, 
$\alpha_A^{-1}\colon \alpha_A(A) \rightarrow A$ is continuous.

(a) By Proposition~\ref{prop:pol:basic}(c), there is $V \in 
\mathcal{N}(\hat A)$ such that $E \subseteq V^\vartriangleright$. By 
Proposition~\ref{prop:ev:cont}(a), $V^\vartriangleright$ is compact in 
$\hat{\hat A}$, and so $V^\vartriangleright \cap \alpha_A(A)$ is 
precompact. Thus, the continuous homomorphic image 
$\alpha^{-1}_A(V^\vartriangleright \cap \alpha_A(A))$ is also
precompact, and contains $\alpha^{-1}_A(E)$. Therefore, 
$\alpha^{-1}_A(E)$ is precompact.

(b) The set $K^\vartriangleright$ is open in $\hat A$, and thus
$K^{\vartriangleright\vartriangleright}$ is compact and equicontinuous in 
$\hat{\hat A}$ (cf. Propositions~\ref{prop:ev:cont}(a) 
and~\ref{prop:pol:basic}(c)). Hence, by (a), 
$K^{\vartriangleright\vartriangleleft}=
\alpha_A^{-1}(K^{\vartriangleright\vartriangleright})$ is precompact in
$A$.

(c) By Propositions~\ref{cor:ev:conthat}, in this case,
$\{K^{\vartriangleright\vartriangleright}\mid K\in\mathcal{K}(A)\}$ is a 
cobase for $\hat{\hat A}$, and so the statement follows from (b).
\end{proof}

\subsection{Quasi-convex compactness}
It is a well-known that the closed convex hull of a weakly compact 
subset of a Banach space is weakly compact. In fact, the property of 
preservation of compactness under formation of closed convex hull 
characterizes completeness in the category of metrizable locally convex 
spaces (cf.~\cite[2.4]{OstWil}). Motivated by this, one says that 
$A \in \mathsf{Ab(Top)}$ has the {\em quasi-convex compactness property} 
(or briefly, {\em $A$ is QCP}) if for every compact subset $K\subseteq A$,
the quasi-convex hull $K^{\vartriangleright\vartriangleleft}$ (of $K$ in 
$A$) is compact.

\begin{proposition}
Let $A \in \mathsf{Ab(Top)}$ be such that $\alpha_A$ is continuous. Then 
$\hat A$ has the quasi-convex compactness property.
\end{proposition}

\begin{proof}
Let $\Phi \subseteq \hat A$ be compact. By Proposition~\ref{prop:ev:cont}(b), 
there is $U \in \mathcal{N}(A)$ such that $\Phi \subseteq 
U^\vartriangleright$, and so $\Phi^{\vartriangleright\vartriangleleft}
\subseteq U^{\vartriangleright\vartriangleright\vartriangleleft}= 
 U^{\vartriangleright}$. Therefore, the closed set 
$\Phi^{\vartriangleright\vartriangleleft}$ is a subset of the compact set
$U^{\vartriangleright}$ (cf. Proposition~\ref{prop:ev:cont}(a)). Hence, 
the result follows.
\end{proof}

\begin{proposition}
Let $A \in \mathsf{Ab(Haus)}$ be a locally quasi-convex group. If

\begin{list}{{\rm (\alph{enumi})}}
{\usecounter{enumi}\setlength{\labelwidth}{25pt}\setlength{\topsep}{-8pt} 
\setlength{\itemsep}{-5pt} \setlength{\leftmargin}{25pt}}
 
\item
$A$ is complete, or

\item
$\alpha_A$ is surjective,

\end{list}
\vspace{6pt}
then $A$ has the quasi-convex compactness property.
\end{proposition}

\begin{proof}
By Proposition~\ref{prop:ev:emb}, $\alpha_A^{-1}\colon \alpha_A(A) 
\rightarrow A$ is continuous.

(a) By Proposition~\ref{prop:lqc:precom}, 
$K^{\vartriangleright\vartriangleleft}$ is precompact in $A$, and it, 
being a closed subspace, is complete. Therefore, by
Theorem~\ref{thm:precomp:cmpct}, it is compact.

(b) Since $K^\vartriangleright$ is open in $\hat A$, 
$K^{\vartriangleright\vartriangleright}$ is compact in $\hat{\hat A}$, and 
therefore $K^{\vartriangleright\vartriangleleft}=
\alpha_A^{-1}(K^{\vartriangleright\vartriangleright})$ is compact, because
$\alpha_A^{-1}$ is continuous.
\end{proof}

\begin{proposition} \label{prop:qcp:alpha}
Let $A \in \mathsf{Ab(Haus)}$ be such that $\alpha_A$ and 
$\alpha_{\hat A}$ are continuous. The following statements are equivalent:

\begin{list}{{\rm (\roman{enumi})}}
{\usecounter{enumi}\setlength{\labelwidth}{25pt}\setlength{\topsep}{-8pt} 
\setlength{\itemsep}{-5pt} \setlength{\leftmargin}{25pt}}

\item
$\alpha_A^{-1}(F)$ is compact for every compact subset $F\subseteq 
\hat{\hat A}$;

\item
$A$ has the quasi-convex compactness property.

\end{list}
\end{proposition}

\begin{proof}
(i) $\Rightarrow$ (ii) is immediate, because 
$K^{\vartriangleright\vartriangleright}$ is compact in $\hat{\hat A}$ and
$K^{\vartriangleright\vartriangleleft}=
\alpha_A^{-1}(K^{\vartriangleright\vartriangleright})$.

(ii) $\Rightarrow$ (i): Let $F \subseteq \hat{\hat A}$ be compact.
Since $F$ is closed and $\alpha_A$ is continuous, $\alpha_A^{-1}(F)$ is 
also closed in $A$. By Corollary~\ref{cor:ev:conthat},
there is a compact subset $K \subseteq A$ such that
$F \subseteq K^{\vartriangleright\vartriangleright}$
(because $\alpha_{\hat A}$ is continuous). Thus,
$\alpha^{-1}(F)\subseteq \alpha_A^{-1}(K^{\vartriangleright\vartriangleright})
=K^{\vartriangleright\vartriangleleft}$, and 
$K^{\vartriangleright\vartriangleleft}$ is compact because $A$ is QCP.
Therefore, $\alpha^{-1}(F)$, being the closed subset of a compact set, is 
compact.
\end{proof}

\section{Special groups and subgroups}


\subsection{Metrizable groups} Every Hausdorff topological group is
completely regular ({\bf nice} reference???), so by Tychonoff's
metrization theorem, if a Hausdorff group is second countable, then it is
metrizable. But for topological groups, second countability is more than
necessary in order to warrant metrizability. Theorem~\ref{thm:Icount:met}
below (originally proved by Kakutani and Birkhoff in 1936) can also be
obtained from well-known results on uniform spaces (cf.~\cite[8.1.10,
8.1.21]{Engel6}).

\begin{ftheorem} \label{thm:Icount:met}
Every first countable Hausdorff topological group is metrizable.
\qed
\end{ftheorem}

We turn to abelian metrizable groups. Recall that a Hausdorff space $X$ is 
a {\em $k$-space} if $F\subseteq X$ is closed in $X$ whenever $F\cap K$ 
is closed for every $K \in \mathcal{K}(X)$. (For details, see 
Appendix~\ref{app:kspace}.)

\begin{ftheorem}[\pcite{Theorem~1}{Chasco3}] \label{thm:met:ksp}
Let $A \in \mathsf{Ab(Met)}$. Then $\hat A$ is a $k$-space.
\end{ftheorem}

Although Theorem~\ref{thm:met:ksp} has a non-commutative generalization 
(cf. \cite{GL3}, \cite[3.4]{GLPHD}), we provide here the original proof 
by Chasco, enriched with a few words of explanation.

\begin{proof}
In order to show that $\hat A$ is a $k$-space, 
let $\Phi \subseteq \hat A$ be such that $\Phi\cap\Xi$ is closed in 
$\Xi$ for every compact subset $\Xi\subseteq \hat A$. We show that 
for every $\zeta \not \in \Phi$, there exists a compact subset 
$K \subseteq A$ such that $(K^\vartriangleright+ \zeta)\cap 
\Phi=\emptyset$, and so $\Phi$ is closed in $\hat A$. 
Without loss of generality, we may assume that $\zeta = 0$.  

Since $A$ is metrizable, $0$ has a decreasing countable base $\{U_n\}$, 
and we set $U_0=A$. Our aim is to construct inductively a family 
$\{F_n\}_{n=0}^\infty$ of finite subsets of $A$ such that for every 
$n\in\mathbb{N}$,
\begin{align}
& F_n \subseteq U_n \mbox{,}  \label{cond:1} \\
& \bigcap\limits_{k=0}^n F_k^\vartriangleright\cap
U_{n+1}^\vartriangleright \cap \Phi =  
\emptyset \mbox{.} \label{cond:2}
\end{align}
By Proposition~\ref{prop:ev:cont}(a), $U_1^\vartriangleright$ is compact 
in $\hat A$, and thus by our assumption, $U_1^\vartriangleright \cap\Phi$
is closed in $U_1^\vartriangleright$. Thus, 
$U_1^\vartriangleright \cap\Phi$ is compact in $\hat A$, and so it
is compact in the pointwise topology carried by $\hom(A,\mathbb{T})$
(which is coarser than the compact-open one). In particular,
$U_1^\vartriangleright \cap\Phi$ is closed in the pointwise topology.
A basic neighborhood of $0$ in the pointwise topology has the form 
$F^\vartriangleright$, where $F \subseteq A$ is finite, so
$0 \not\in U_1^\vartriangleright \cap\Phi$ implies that there exists
a finite subset $F_0 \subseteq A=U_0$  such that 
$F_0^\vartriangleright \cap U_1^\vartriangleright \cap\Phi=\emptyset$.
This completes the proof for $n=0$.

Suppose that we have already constructed $F_0,\ldots,F_{n-1}$ such that
(\ref{cond:1}) and (\ref{cond:2}) hold. For each $x \in U_n$, put
\begin{align}
\Delta_x & = \bigcap\limits_{k=0}^{n-1}
F_k^\vartriangleright \cap \{x\}^\vartriangleright \cap 
U_{n+1}^\vartriangleright \cap \Phi \mbox{.} \\
\intertext{Since $U_n \supseteq U_{n+1}$, one has 
$U_n^\vartriangleright \subseteq U_{n+1}^\vartriangleright$, and thus, by 
the inductive hypothesis,}
\bigcap\limits_{x \in U_n} \Delta_x  &=
\bigcap\limits_{k=0}^{n-1}
F_k^\vartriangleright \cap U_n^\vartriangleright \cap
U_{n+1}^\vartriangleright \cap \Phi =
\bigcap\limits_{k=0}^{n-1}
F_k^\vartriangleright \cap U_n^\vartriangleright \cap \Phi =
\emptyset
\mbox{.} \\
\intertext{The $\Delta_x$ are closed subsets of the compact space 
$U_{n+1}^\vartriangleright$, and their intersection is empty. Therefore, 
there must be a finite subset $F_n \subseteq U_n$ such that 
$\bigcap\limits_{x \in F_n} \Delta_x = \emptyset$. Hence,}
\bigcap\limits_{x \in F_n} \Delta_x & =  
\bigcap\limits_{k=0}^{n-1}
F_k^\vartriangleright \cap F_n^\vartriangleright \cap
U_{n+1}^\vartriangleright \cap \Phi =
\emptyset,
\end{align}
as desired.
Set  $K = \bigcup\limits_{n=0}^\infty F_n \cup \{0\}$. It follows from
(\ref{cond:1}) that $K\backslash U_n$ is finite for every 
$n \in  \mathbb{N}$, and thus $K$ is sequentially compact (because every 
sequence without a constant subsequence must converge to $0 \in K$). 
Therefore, $K$ is compact, because $A$ is metrizable. By the construction 
of $K$, $K^\vartriangleright \cap U_n^\vartriangleright \cap 
\Phi=\emptyset$ for every $n \in \mathbb{N}$. By 
Proposition~\ref{prop:pol:basic}(a), 
$\hat A= \bigcup\limits_{n=1}^\infty U_n^\vartriangleright$, and therefore
$K^\vartriangleright \cap \Phi=\emptyset$, as desired. 
\end{proof}

\begin{corollary} \label{cor:met:alpha}
Let $A \in \mathsf{Ab(Met)}$. Then:

\begin{list}{{\rm (\alph{enumi})}}
{\usecounter{enumi}\setlength{\labelwidth}{25pt}\setlength{\topsep}{-8pt} 
\setlength{\itemsep}{-5pt} \setlength{\leftmargin}{25pt}}

\item
$\alpha_A$ is continuous;

\item
$\alpha_{\hat A}$ is continuous;

\item
{\rm (\cite[Corollary~1]{Chasco3})}
$\hat{\hat A}$ is complete and metrizable.

\end{list}
\end{corollary}

\begin{proof}
(a) Since $A$ is metrizable, its topology is determined by convergent 
sequences $x_n \longrightarrow x_0$. For such a sequence, 
$\{x_n\mid n \in \mathbb{N}\} \cup \{x_0\}$ is compact, and thus $A$  is 
a $k$-space. Therefore, by Theorem~\ref{thm:ev:kcont}(b), $\alpha_A$ is 
continuous.

(b) By Theorem~\ref{thm:met:ksp}, $\hat A$ is a $k$-space,
and so, by Theorem~\ref{thm:ev:kcont}(b), $\alpha_{\hat A}$ 
is continuous.

(c) By Proposition~\ref{prop:ev:cont}(b), 
$\{U^{\vartriangleright\vartriangleright} \mid U \in \mathcal{N}(A)\}$
is a base for $\hat{\hat A}$ at $0$, because $\alpha_A$ is continuous.
Thus, $\hat{\hat A}$ is first-countable, because $A$ is so. (Recall
that for topological groups, being first-countable is equivalent to 
being metrizable.) By Theorem~\ref{thm:met:ksp}, $\hat A$ is a $k$-space, and 
thus $\mathscr{C}(\hat A,\mathbb{T})$
is complete, because $\mathbb{T}$ is so (cf.~\cite[7.12]{Kelley}).
Since $\hat{\hat A}$ is a closed subspace of 
$\mathscr{C}(\hat A,\mathbb{T})$, this completes the proof.
\end{proof}

A combination of Proposition~\ref{prop:qcp:alpha} and 
Corollary~\ref{cor:met:alpha} yields:

\begin{corollary} \label{cor:met:qcp}
Let $A \in \mathsf{Ab(Met)}$. Then $A$ has the quasi-convex compactness 
property if and only if $\alpha_A^{-1}(F)$ is compact in $A$ for every 
compact $F \subseteq \hat{\hat A}$.
\qed
\end{corollary}

\begin{fcorollary}[\pcite{Theorem~2}{Chasco3}, \pcite{3.7}{GLPHD}] 
\label{cor:met:dense}
Let $A \in \mathsf{Ab(Met)}$, and let $D \leq A$ be its dense subgroup.
Then $\hat D=\hat A$ as topological groups.
\end{fcorollary}

\begin{proof}
By Theorem~\ref{thm:dense:sbgrp}, $\hat D$ and $\hat A$ have the same 
underlying groups and compact subsets. On the other hand, by 
Theorem~\ref{thm:met:ksp}, both $\hat D$ and $\hat A$ are $k$-spaces,
because $D$ and $A$ are metrizable. Therefore, $\hat A$ and $\hat D$ have 
the same topology, because the topology of a $k$-space is determined by its 
compact subsets. 
\end{proof}

\begin{theorem}
Let $A \in \mathsf{Ab(Met)}$. The following statements are equivalent:

\begin{list}{{\rm (\roman{enumi})}}
{\usecounter{enumi}\setlength{\labelwidth}{25pt}\setlength{\topsep}{-8pt} 
\setlength{\itemsep}{-5pt} \setlength{\leftmargin}{25pt}}

\item
$A$ is locally quasi-complete and complete;

\item
$\alpha_A$ is a closed embedding;


\item
$A$ has the quasi-convex compactness property and $\alpha_A$ is injective.

\end{list}
\end{theorem}

Ostling and Wilansky showed that a locally convex metrizable vector space 
is complete if and only if the absolutely convex closure of 
compact subset is compact (cf.~\cite[2.4]{OstWil}).
Hern\'andez proved a far reaching generalization of this 
result for metrizable groups, namely, that a metrizable LQC group 
is complete if and only if it has QCP (cf.~\cite[Theorem~2]{Hern2}). 
The same
result also appears in a paper by Bruguera and Mart\'\i n-Peinador,
who used Corollary~\ref{cor:met:dense} in order to simplify the proof
(cf.~\cite[9]{BrugMar2}). It escaped the attention of these
authors that if $\alpha_A$ is injective, then QCP is sufficient
to warrant completeness, so the first half of the proof of 
(iii) $\Rightarrow$ (i) appears to be new, while its second half
uses the ideas of Bruguera and Mart\'\i n-Peinador.

\begin{proof}
(i) $\Rightarrow$ (ii): Since $A$ is LQC (and metrizable, so Hausdorff), 
$\alpha_A$ is an embedding by Proposition~\ref{prop:ev:emb}. Furthermore, 
if $S \subseteq A$ is closed in $A$, then $S$ is complete (because $A$ is 
complete), and thus $\alpha_A(S)$ is complete. Therefore, $\alpha_A(S)$ is 
closed in $\hat{\hat A}$, as desired.

(ii) $\Rightarrow$ (iii): Clearly, $\alpha_A$ is injective. Let 
$K \subseteq A$ be compact. Then $K^{\vartriangleright\vartriangleright}$ 
is compact in $\hat{\hat A}$, and since $\alpha_A(A)$ is closed,
$K^{\vartriangleright\vartriangleright} \cap \alpha_A(A)$ is also compact.
Therefore, $K^{\vartriangleright\vartriangleleft}=
\alpha_A^{-1}(K^{\vartriangleright\vartriangleright} \cap \alpha_A(A))$ is 
compact too.

(iii) $\Rightarrow$ (i): Since $\alpha_A$ is injective, 
$\alpha_A^{-1}\colon \alpha_A(A)\rightarrow A$ is well-defined. By 
Corollary~\ref{cor:met:qcp}, QCP implies that $\alpha_A^{-1}(F)$ is compact 
for every compact subset $F \subseteq \alpha_A(A)$. So, if
$K \subseteq \alpha_A^{-1}(F)$ is closed, then $K$ is also compact, and 
thus $(\alpha_A^{-1})^{-1}(K)=\alpha_A(K)$ is a compact subset of $F$, 
because $\alpha_A$ is continuous. Therefore, 
$(\alpha_A^{-1})^{-1}(K)$ is closed, and hence, $\alpha_A^{-1}$ 
is $k$-continuous. By Corollary~\ref{cor:met:alpha}(c), $\hat{\hat A}$ is 
metrizable, and so is its subspace $\alpha_A(A)$. In particular,
$\alpha_A(A)$ is a $k$-space, and therefore $\alpha_A^{-1}$ is 
continuous. This shows that $\alpha_A$ is an embedding, and hence,
by Proposition~\ref{prop:ev:emb}(a), $A$ is LQC.

In order to show that $A$ is complete, let $B$ be the completion of $A$.
By Corollary~\ref{cor:met:dense}, $\hat A= \hat B$ (as topological 
groups), because $A$ is a dense subgroup of the metrizable group $B$.
Let $\{a_n\}$ be a Cauchy-sequence in $A$. Then $a_n 
\longrightarrow b$ for some $b \in B$, and 
$K_1=\{a_n\mid n \in \mathbb{N}\}\cup \{b\}$ is compact in $B$. So,
$K_1^\vartriangleright \in \mathcal{N}(\hat B)=\mathcal{N}(\hat A)$, and 
by Proposition~\ref{prop:pol:basic}(b), there is a compact subset 
$K \subseteq A$ such that $K^\vartriangleright \subseteq 
K_1^\vartriangleright$. Thus, 
$\{a_n\mid n \in \mathbb{N}\} \subseteq 
K_1^{\vartriangleright\vartriangleleft}\subseteq 
K^{\vartriangleright\vartriangleleft}$ (where $(-)^\vartriangleleft$ is 
taken with respect to $A$). Since $A$ has QCP, 
$K^{\vartriangleright\vartriangleleft}$ is compact, and so 
sequentially compact (because $A$ is metrizable). Therefore, $\{a_n\}$ has 
a convergent subsequence in $K^{\vartriangleright\vartriangleleft}$, which 
shows that $b \in K^{\vartriangleright\vartriangleleft}\subseteq A$.
Hence, $A$ is complete, as desired.
\end{proof}


\subsection{Compact and discrete groups} Let $A \in \mathsf{Ab(Top)}$.
If $A$ is compact, then by Proposition~\ref{prop:pol:basic}(b), 
$\{0\}=A^\vartriangleright\in\mathcal{N}(\hat A)$. On the 
other hand, if $A$ is discrete, then $\{0\}\in \mathcal{N}(A)$, and so
by Proposition~\ref{prop:ev:cont}(a), $\hat A = \{0\}^\vartriangleright$ 
is compact. Repeating these arguments for $\hat A$ yields:

\begin{samepage}

\begin{flemma} \label{lemma:spec:cd}
Let $A \in \mathsf{Ab(Top)}$.

\begin{list}{{\rm (\alph{enumi})}} 
{\usecounter{enumi}\setlength{\labelwidth}{25pt}\setlength{\topsep}{-8pt} 
\setlength{\itemsep}{-5pt} \setlength{\leftmargin}{25pt}}
 
\item
If $A$ is compact, then $\hat A$ is discrete.

\item
If $A$ is discrete, then $\hat A$ is compact.

\item
If $A$ is compact, then so is $\hat{\hat A}$.

\item
If $A$ is discrete, then so is $\hat{\hat A}$. 
\qed

\end{list}
\end{flemma}

\end{samepage}

\begin{ftheorem} \label{thm:P:cd}
Let $A \in \mathsf{Ab(Top)}$. If $A$ is 

\begin{list}{{\rm (\alph{enumi})}}
{\usecounter{enumi}\setlength{\labelwidth}{25pt}\setlength{\topsep}{-8pt} 
\setlength{\itemsep}{-5pt} \setlength{\leftmargin}{25pt}}

\item discrete, or 

\item compact Hausdorff, 

\end{list}
\vspace{6pt}
then $\alpha_A$ is an isomorphism of topological groups.
\end{ftheorem}

\begin{proof}
(a) Let $a \in A$ be a non-zero element, and 
let $r \in \mathbb{T}$ be such that $o(a)=o(r)$. This defines a 
group homomorphism $\chi_a\colon \langle a \rangle \rightarrow 
\mathbb{T}$ such that $\chi_a(a)\neq 0$.
Since $\mathbb{T}$ is an injective abelian group, $\chi_a$ 
extends to $\bar\chi_a\colon A \rightarrow \mathbb{T}$, 
and $\bar \chi_a \in \hat A$, because $A$ is discrete. Thus,
$\alpha_a(\bar\chi_a)=\bar\chi_a(a)\neq 0$, which
shows that
$\alpha_A$ is injective. By Lemma~\ref{lemma:spec:cd}(d), $\hat{\hat A}$ 
is discrete, and thus $\alpha_A$ is an embedding.

Suppose that $A$ is finitely generated. Then it decomposes into the 
direct sum of cyclic groups $A=\langle a_1 \rangle \oplus \cdots \oplus 
\langle a_n \rangle$. Finite cyclic groups are self-dual, and 
$\hat{\mathbb{Z}}=\mathbb{T}$ and $\hat{\mathbb{T}}=\mathbb{Z}$. Thus,  
surjectivity of $\alpha_A$ follows, because the Pontryagin dual is an 
additive functor. (More details???)

In the general case, let $\xi \in \hat{\hat A}$. By 
Proposition~\ref{prop:pol:basic}(a), there is
$W\in \mathcal{N}(\hat A)$  such that 
$\xi(W)\subseteq\nolinebreak \Lambda_1$, and $W$ can be chosen to have 
the form $F^\vartriangleright$ where $F \subseteq A$ is finite. Put
$B=\langle F \rangle$.  By Lemma~\ref{lemma:sub:ann},
$B^\vartriangleright$ is a subgroup of $\hat A$, and 
$B^{\vartriangleright\vartriangleright}\cong 
\widehat{\hat A/B^\vartriangleright}$. Since $\mathbb{T}$ is an injective 
abelian group, every homomorphism $\chi\colon B \rightarrow \mathbb{T}$ 
extends to a homomorphism $\bar\chi\colon A \rightarrow \mathbb{T}$, and 
thus the dual $\hat\iota_B\colon \hat A \rightarrow \hat B$ of the 
inclusion $\iota_B\colon B \rightarrow A$ is surjective. The map 
$\hat\iota_B$ is also closed, because $\hat A$ is compact. Therefore, 
$\hat B \cong \nolinebreak \hat A /B^\vartriangleright$, and hence
$B^{\vartriangleright\vartriangleright}\cong
\widehat{\hat A/B^\vartriangleright} \cong \hat{\hat B}$.
Since $F \subseteq B$, one has 
$\xi(B^\vartriangleright) \subseteq \xi(F^\vartriangleright) \subseteq 
\Lambda_1$, and so $\xi \in \nolinebreak
B^{\vartriangleright\vartriangleright}\cong \nolinebreak
\hat{\hat B}$. The group $B$ 
is finitely generated, and thus by what we have shown so far, there is 
$b \in B$ such that $\xi=\alpha_B(b)$. Therefore, $\xi=\alpha_A(b)$.

(b) Since $A$ is compact Hausdorff, 
it is  a $k$-space, and so $\alpha_A$ is continuous (cf. 
Theorem~\ref{thm:ev:kcont}). By the Peter-Weyl theorem, for every
non-zero $a \in A$, there is $\chi\in \hat A$ such that 
$\chi(a)\neq 0$ (cf. \cite[Thm.~32]{Pontr}). Thus, $\alpha_A$ is 
injective, and therefore it is a closed embedding (because $A$
is compact). In order to show that 
$\alpha_A$ is surjective, assume the contrary, that is, assume
that $\alpha_A(A)$ is a proper closed subgroup of $\hat{\hat A}$. By 
Lemma~\ref{lemma:spec:cd}(c), $\hat{\hat A}$ is compact, and thus so is
the quotient $\hat{\hat A}/\alpha_A(A)$. By the Peter-Weyl theorem 
applied to $\hat{\hat A}/\alpha_A(A)$, there exists
a non-zero continuous character 
$\psi\colon\hat{\hat A}/\alpha_A(A)\rightarrow\mathbb{T}$, which 
induces $\zeta\in \hat{\hat{\hat A}}$ such that 
$\zeta(\alpha_A(A))=\{0\}$ and $\zeta \neq 0$. Since $\hat A$ is discrete 
(cf.~Lemma~\ref{lemma:spec:cd}(a)), by (a), $\alpha_{\hat A}$ is an 
isomorphism. In particular, there is a non-zero $\chi\in\hat A$ such  that 
$\alpha_{\hat A}(\chi)=\zeta$. Therefore, for $a \in A$,
\begin{align}
\chi(a)=(\alpha_A(a))(\chi)=(\alpha_{\hat A}(\chi))(\alpha_A(a))=
\zeta (\alpha_A(a)) \in \zeta(\alpha_A(A))=\{0\}.
\end{align}
Hence, $\chi=0$ (and so $\zeta=0$) contrary to our assumption, which 
completes the proof.
\end{proof}

\subsection{Bohr-compactification} The \v Cech-compactification $\beta X$
of a Tychonoff space $X$ has the property that every continuous 
map $f\colon X \rightarrow K$ into a compact Hausdorff space $X$ extends 
to $\beta X$, and thus factors uniquely through the dense embedding $X 
\rightarrow \beta X$. It appears to be less known, however, that $\beta X$ 
exists even when $X$ is not Tychonoff, in which case the map $X 
\rightarrow \beta X$ is only continuous but need not be injective or open 
onto its image. (Every continuous map of $X$ into a compact Hausdorff 
space still factors uniquely through $\beta X$.) In a categorical 
language, one says that the category $\mathsf{HComp}$ of compact Hausdorff 
spaces and their continuous maps is a reflective subcategory of 
$\mathsf{Top}$.

Similarly to the relationship between $\mathsf{Top}$ and $\mathsf{HComp}$, 
the full subcategory $\mathsf{Grp(HComp)}$ of 
compact Hausdorff groups is reflective in $\mathsf{Grp(Top)}$ (topological 
groups and their continuous homomorphisms). The reflection is called the 
{\em Bohr-compactification}, and is denoted by 
$\rho_G\hspace{-0.2pt}\colon \hspace{-0.4pt} G \rightarrow 
\nolinebreak bG$. We show its existence for abelian groups, where the 
construction is rather simple

\begin{theorem}
For every $A \in \mathsf{Ab(Top)}$ there is $bA \in \mathsf{Ab(HComp)}$ 
and a continuous homomorphism $\rho_A\colon A \rightarrow bA$ such that 
every continuous homomorphism $\varphi\colon A\rightarrow K$ 
into $K \in\nolinebreak \mathsf{Ab(HComp)}$ factors uniquely through 
$\rho_A$:
\begin{equation}
\bfig
\Vtriangle(0,0)/->`->`<--/<300,350>%
[A`K`bA;\varphi`\rho_A`\exists!\tilde\varphi]
\efig
\end{equation}
Moreover, $\widehat{bA}=(\hat A)_d$, and $\ker \rho_A = \ker \alpha_A$.
\end{theorem}

\begin{proof}
Since $bA$ is supposed to be compact, by 
Lemma~\ref{lemma:spec:cd}(a), $\widehat{bA}$ must be discrete. 
Applying the universal property of $bA$ to $K=\mathbb{T}$ and 
$\varphi=\chi \in \hat A$ yields that $A$ and $bA$ have the same 
continuous characters, that is, $\widehat{bA}=(\hat A)_d$
(the subscript $d$ stands for the discrete topology). Using this 
necessary condition as a definition, one set $bA=\widehat{(\hat A)_d}$
and $(\rho_A(a))(\chi)=\chi(a)$ for every $a \in A$ and $\chi\in 
(\hat A)_d$.
By Proposition~\ref{prop:pol:basic}(b), a basic open set in $bA$ has the 
form $\Phi^\vartriangleright$, where $\Phi \subseteq \hat A$ is finite. 
Thus, $\rho_A^{-1}(\Phi^\vartriangleright) = \bigcap\limits_{\chi\in\Phi}
\chi^{-1}(\Lambda_1)$ is a neighborhood of $0$ in $A$, and 
therefore $\rho_A$ is continuous. If $\varphi\colon A \rightarrow K$ is a 
continuous homomorphism, then it induces $\hat\varphi\colon \hat K 
\rightarrow \hat A$. Since $\hat K$ is discrete (cf. 
Lemma~\ref{lemma:spec:cd}(a)), 
$(\hat\varphi)_d\colon \hat K \rightarrow (\hat A)_d$ is also continuous, 
and so $\widehat{(\hat\varphi)_d}\colon \widehat{(\hat A)_d} \rightarrow
\hat{\hat K}$ is continuous as well. By Theorem~\ref{thm:P:cd},
$\hat{\hat K}\cong K$, and therefore $\tilde\varphi=
\widehat{(\hat\varphi)_d}$. This consideration shows the uniqueness 
of $\tilde \varphi$ too, because $\tilde \varphi$ is uniquely determined by 
its dual. Since $\alpha_A(a)$ and $\rho_A(a)$ are both evaluations at 
$a$, the last statement follows.
\end{proof}

A subgroup $H \leq A$ such that $H^{\vartriangleright\vartriangleleft}=H$ 
is said to be {\em dually closed} in $A$; if every continuous character 
of $H$ extends continuously to $A$, we say that $H$ is {\em dually 
embedded} in $A$. We saw in Theorem~\ref{thm:dense:sbgrp}(a) that every 
dense subgroup of a topological group is dually embedded.

\begin{corollary} \label{cor:comp:dually}
Let $A \in \mathsf{Ab(Top)}$ be such that $\alpha_A$ is injective, and let 
$K$ be a compact subgroup of $A$. Then $K$ is dually closed and dually 
embedded in $A$.
\end{corollary}

\begin{proof}
Since $\alpha_A$ is injective, $A$ must be Hausdorff. First, suppose that 
$A$ is compact. In order to show that $K$ is dually closed, let 
$x \in A \backslash K$. The quotient $A/K$ is also compact, and 
the image $\bar x$ of $x$ in $A/K$ is non-zero. Thus, by 
Theorem~\ref{thm:P:cd}, there is $\bar\chi \in \widehat{A/K}$ such that 
$\bar\chi(\bar x) \neq 0$. The character $\bar \chi$ induces
$\chi \in \hat A$ by setting $\chi(a)=\bar\chi(\bar a)$ (where $\bar a$ 
stands for the image of $a$ in $A/K$), and $K \subseteq \ker \chi$, 
while $\chi(x)\neq 0$. Therefore, $\chi \in K^\vartriangleright$, and 
hence $x \not\in K^{\vartriangleright\vartriangleleft}$. (Note that by 
Lemma~\ref{lemma:sub:ann}, $K^\vartriangleright$ is a subgroup of $\hat A$, 
and so $K^{\vartriangleright\vartriangleleft}$ coincides with its 
annihilator in $A$.) This shows that 
$K=K^{\vartriangleright\vartriangleleft}$. In order to show that $K$ is 
dually embedded, consider the discrete group $\hat A$. 
By  Lemma~\ref{lemma:sub:ann}(b), 
$\hat\pi_K\colon \widehat{\hat A/ K^\vartriangleright} \rightarrow 
\hat{\hat A}$ is injective and its image is 
$K^{\vartriangleright\vartriangleright}$. Since $\hat{\hat A}\cong A$
(cf.~Theorem~\ref{thm:P:cd}), $K^{\vartriangleright\vartriangleright}\cong
K^{\vartriangleright\vartriangleleft}=K$. The quotient 
$\hat A/ K^\vartriangleright$ is discrete, and so
$\widehat{\hat A/ K^\vartriangleright}$ is compact. 
Therefore, $\hat\pi_K\colon \widehat{\hat A/ K^\vartriangleright} 
\rightarrow A$ is an embedding onto $K$, that is,
$K \cong \widehat{\hat A/ K^\vartriangleright}$. Hence, applying 
Theorem~\ref{thm:P:cd} to the compact group $\hat A/ K^\vartriangleright$ 
yields $\hat K \cong \widehat{\widehat{\hat A/ K^\vartriangleright}} \cong
\hat A/ K^\vartriangleright$. In particular, 
$\hat \iota_A \colon \hat A \rightarrow \nolinebreak\hat K$ given by 
restriction to $K$ is surjective, as desired.

In the general case, since $\rho_A$ continuous, 
$\rho_A(K)$ is a compact subgroup of $bA$. Therefore, by what we have
shown so far, it is dually closed and dually embedded in $bA$. Thus,
for every $x \in bA\backslash \rho_A(K)$, there is $\chi \in 
\widehat{bA}=(\hat A)_d$ such that $\chi(\rho_A(K))=0$ and 
$\chi(x)\neq 0$. In particular, this holds for every $x=\rho_A(a)$, 
where $a \in A \backslash K$. (Since $\rho_A$ is injective, $a \in A 
\backslash K$ implies $x \in bA\backslash \rho_A(K)$.) Therefore,
$\chi \in K^\vartriangleright$, and so $a \not\in 
K^{\vartriangleright\vartriangleleft}$. Hence, 
$K=K^{\vartriangleright\vartriangleleft}$, as desired. For the second 
statement, observe that $K \cong \rho_A(K)$ (because $\rho_A$ is 
injective, $K$ is compact, and $bA$ is Hausdorff). Thus, if
$\chi \in \hat K=\widehat{\rho_A(K)}$, then it admits a continuous 
extension $\psi$ to $bA$, because $\rho_A(K)$ is dually embedded in 
$bA$. Therefore, $\psi\circ \rho_A \in \hat A$ is a continuous extension 
of $\chi$. Hence, $K$ is dually embedded in $A$, as desired.
\end{proof}

\subsection{Open and compact subgroups} Open subgroups of abelian groups 
are closely related to compact subgroups of the Pontryagin dual, 
and vice versa. We turn to 
establishing the properties of open and compact subgroup of abelian 
groups. Our presentation of these results was strongly influenced by the 
work of Banaszczyk, Chasco, Mart\'\i n-Peinador, and Bruguera (cf. 
\cite{Chasco1} and \cite{BrugMar}).

\begin{fproposition}[\pcite{2.2}{Chasco1}]  \label{prop:open:subgp}
Let $A\in\mathsf{Ab(Top)}$, $U\leq A$ be an open subgroup, and consider 
the exact sequence
\begin{align}
\hspace{28pt} 0 \to U \to^{\iota_U} A \to^{\pi_U} A/U \to 0.
\end{align}

\begin{list}{{\rm (\alph{enumi})}}
{\usecounter{enumi}\setlength{\labelwidth}{25pt}\setlength{\topsep}{-8pt} 
\setlength{\itemsep}{-5pt} \setlength{\leftmargin}{25pt}}


\item 
The map $\hat\pi_U\colon \widehat{A/U} \rightarrow \hat A$ is an 
embedding and $\widehat{A/U} \cong U^\vartriangleright$.

\item
The subgroup $U$ is dually embedded in $A$.
In other words, $\hat\iota_U\colon \hat A \rightarrow \hat U$ is 
surjective.

\item
$\hat\iota_U$ is open, and thus a quotient.

\item
The induced sequence
\begin{align}
0 \to \widehat{A/U} \to^{\hat \pi_U} \hat A \to^{\hat\iota_U}   
\hat U \to 0
\end{align}
is exact.

\end{list}
\end{fproposition}

\begin{proof}
Note that $A/U$ is discrete, because $U$ is open, and so
$\widehat{A/U}$ is compact by Lemma~\ref{lemma:spec:cd}(b).

(a) By Lemma~\ref{lemma:sub:ann}(b), $\hat\pi_U$ is injective, and its 
image is $U^\vartriangleright$. Thus, $\hat\pi_U$ is an embedding, because
it is continuous, its domain is compact and its codomain is Hausdorff. 

(b) Let $\chi \in \hat U$, and let $V\in \mathcal{N}(U)$ be such that
$\chi(V) \subseteq \Lambda_1$ (cf. 
Proposition~\ref{prop:pol:basic}(a)). Since $\mathbb{T}$ is an injective 
abelian group, $\chi$ admits an extension 
$\psi\colon A \rightarrow \mathbb{T}$. Thus, 
$\psi(V)=\chi(V) \subseteq \Lambda_1$, and 
$V \in \mathcal{N}(A)$, because $U$ is open in $A$. Therefore, 
by Proposition~\ref{prop:pol:basic}(a), $\psi$ is continuous on $A$.

(c) In order to make he proof more transparent, we decomposed it into 
three simple steps.

{\em Step 1:} $A=\langle U, b\rangle$ for some $b \in A$. Let $l$ be the 
order of $b$ in $A/U$ (possibly infinite). Given a compact subset
$K\subseteq A$, it can be covered by finitely many translates
$k_1b+U,\ldots, k_mb+U$  of $U$, and without loss of generality we may 
assume that $0 \leq k_i < l$. Since $U$ is an open subgroup, it is also 
closed, and so  $(K -k_ib) \cap  U$ is compact. Thus,
$C=((K -k_1b)\cup \ldots \cup (K-k_mb)) \cap  U$ is a compact subset of 
$U$, and $K \subseteq (k_1b+ C) \cup \ldots \cup (k_m b + C)$. Without 
loss of generality, we may assume that $lb \in C$ if $l\neq \infty$.
Set $W=\{\chi \in \hat U\mid \chi(C) \subseteq \Lambda_2\}$, and 
we show that $W \subseteq \hat\iota_U(K^\vartriangleright)$. To that end, 
let $\chi \in W$. If $\bar\chi \in \hat A$ is an extension of $\chi$, then
\begin{align} \label{eq:chi:cover}
\bar\chi(K) \subseteq \bar\chi(k_1b+ C) \cup \ldots \cup \bar\chi(k_mb+C)
\subseteq (k_1\bar\chi(b) +\Lambda_2)\cup \ldots 
(k_m\bar\chi(b) +\Lambda_2).
\end{align}
If $l$ is finite, then $lb \in C$, and so $\chi(lb) \in \Lambda_2$.
Let $r \in\Lambda_{2l} \subseteq \mathbb{T}$ be the closest point 
to $0$ such that $lr=\chi(lb)$ (if $l=\infty$, then $r=0$).  
Set $\bar\chi(u+nb)=\chi(u)+nr$. The character $\chi$ is 
continuous on $A$ because $\bar\chi\rest_U=\chi$ is continuous, 
and $U$ is an open subgroup. Furthermore, one has 
$k_i \bar\chi(b) \in \Lambda_2$ for each $i$, because
$0 \leq k_i < l$ and $r \in\Lambda_{2l}$. Therefore, by 
(\ref{eq:chi:cover}),
$\bar\chi(K) \subseteq \Lambda_2 + \Lambda_2 \subseteq \Lambda_1$, which 
means that $\bar\chi \in K^\vartriangleright$. 
Hence, $W \subseteq \hat\iota_U(K^\vartriangleright)$, as desired.

{\em Step 2:} $A/U$ is finitely generated, that is, 
$A=\langle U,b_1,\ldots,b_n\rangle$ for $b_1,\ldots b_n \in A$. Set
$U_0=U$ and $U_k=\langle U_{k-1},b_k\rangle$. Each $U_k$ is an open 
subgroup of $U_{k+1}$ (and of $A$), and thus by step 1, each 
$\hat\iota_{U_k}\colon\hat{U}_{k+1}\rightarrow\nolinebreak\hat{U}_{k}$ is 
open. Therefore, $\hat\iota_U=\hat \iota_{U_{k-1}} \circ \ldots\circ 
\hat\iota_{U_0}$ is open.

{\em Step 3:} In the general case, let $K \subseteq A$ be compact.
Then $K$ can be covered by finitely many translates of $U$, 
$K \subseteq (b_1+U)\cup\cdots\cup (b_n+U)$. Put 
$U^\prime=\langle  U,b_1,\ldots,b_n\rangle$. Since $K \subseteq U^\prime$, 
the image $\hat\iota_{U^\prime}(K^\vartriangleright)$ coincides with 
the polar of $K$ in $\hat{U}^\prime$ with respect to $U^\prime$, and thus 
open. By step 2, $\hat\iota_U\colon \hat{U}^\prime \rightarrow \hat{U}$ is 
open. Therefore, $\hat\iota_U(K^\vartriangleright)$ is open  in $\hat U$, 
as desired.

(d) Exactness at $\widehat{A/U}$ and $\hat U$ were shown in (b) and (c), 
respectively. Exactness at $\hat A$ also follows from (b), because
$\operatorname{Im} \hat\pi_U=U^\vartriangleright=\ker \hat\iota_U$.
\end{proof}

\begin{fproposition} \label{prop:compact:subgr}
Let $A\in\mathsf{Ab(Top)}$, $K\leq A$ be a compact subgroup, and consider the 
exact sequence
\begin{align} \label{eq:compact:subgr}
\hspace{28pt} 0 \to K \to^{\iota_K} A \to^{\pi_K} A/K \to 0.
\end{align}

\begin{list}{{\rm (\alph{enumi})}}
{\usecounter{enumi}\setlength{\labelwidth}{25pt}\setlength{\topsep}{-8pt} 
\setlength{\itemsep}{-5pt} \setlength{\leftmargin}{25pt}}


\item 
The map $\hat\pi_K\colon \widehat{A/K} \rightarrow \hat A$ is an 
embedding and $\widehat{A/K} \cong K^\vartriangleright$.

\item
The map $\hat\iota_K\colon \hat A \rightarrow \hat K$ is open onto its 
image.

\item
The induced sequence
\begin{align}
0 \to \widehat{A/K} \to^{\hat \pi_K} \hat A \to^{\hat\iota_K} \hat K 
\end{align}
is exact.

\item
If $\alpha_A$ is injective, then $\hat\iota_K\colon\hat A\rightarrow\hat K$
is surjective, and so $\hat\iota_K$ is a quotient map.

\end{list}
\end{fproposition}

\begin{proof}
Note that by Lemma~\ref{lemma:spec:cd}(a), $\hat K$ is compact.

(a) A basic neighborhood of $0$ in $\widehat{A/K}$ is of the form 
$L^\vartriangleright$, where $L \subseteq A/K$ is compact (cf. 
Proposition~\ref{prop:pol:basic}(b)). By 
Lemma~\ref{lemma:comp:proper-new}(b), 
$\pi_K^{-1}(L)$ is a compact subset of $A$, and thus
$\pi_K^{-1}(L)^\vartriangleright$ is a basic neighborhood of $0$ in $\hat 
A$. Without loss of generality, we may assume that $0 \in L$, and then by 
Lemma~\ref{lemma:sub:ann}(c), $\pi_K^{-1}(L)^\vartriangleright=
\hat\pi_K(L^\vartriangleright)$, as desired.

(b) Since $\ker \hat\iota_K=K^\vartriangleright$ is open in $\hat A$ and 
$\hat K$ is discrete, the statement is obvious.

(c) Exactness at $\widehat{A/K}$ and at $\hat A$ follows from (b), because 
$\operatorname{Im} \hat\pi_K=K^\vartriangleright=\ker \hat\iota_K$.

(d) By Corollary~\ref{cor:comp:dually}, $K$ is dually embedded in $A$,
which completes the proof.
\end{proof}

\begin{ftheorem}[\pcite{2.3}{Chasco1}] \label{thm:open:subgr}
Let $A \in \mathsf{Ab(Top)}$, and let $U\leq A$ be an open subgroup. 
Then:

\begin{list}{{\rm (\alph{enumi})}}
{\usecounter{enumi}\setlength{\labelwidth}{25pt}\setlength{\topsep}{-8pt} 
\setlength{\itemsep}{-5pt} \setlength{\leftmargin}{25pt}}

\item
$\alpha_U$ is injective (resp., surjective) if and only if $\alpha_A$ 
is injective (resp., surjective);

\item
$\alpha_U$ is an isomorphism of topological groups if and only if 
$\alpha_A$ is so.

\end{list}
\end{ftheorem}

\begin{proof}
By Proposition~\ref{prop:open:subgp}, the exact sequence
\begin{align}
0 \to<285> U \to<285>^{\iota_U} &A \to^{\pi_U} A/U \to 0  \\
\intertext{gives rise to an exact sequence}
0\to\widehat{A/U}\to^{\hat\pi_U}
&\hat A\to<285>^{\hat\iota_U}\hat U\to<285> 0 \\
\intertext{with $\widehat{A/U}$ compact, $\hat\pi_U$ an embedding, and 
$\hat\iota_U$ a quotient map. Thus, the conditions of 
Proposition~\ref{prop:compact:subgr} are fulfilled (by 
Remark~\ref{rem:MAP}(3), $\alpha_{\hat A}$ is injective), and so}
0 \to<285> \hat{\hat{U}} \to<285>^{\hat{\hat\iota}_U} &
\hat{\hat A} \to^{\hat{\hat\pi}_U} \widehat{\widehat{A/U}} \to 0
\end{align}
is also exact, $\hat{\hat\iota}_U$ is an embedding, and $\hat{\hat\pi}_U$ 
is a quotient map. Since the evaluation homomorphism $\alpha$ is a natural 
transformation, we obtain a commutative diagram with exact 
rows:
\begin{align}
\bfig
\morphism(-796,200)<405,0>[0`U;]
\morphism(-391,200)<391,0>[U`A;\iota_U]
\morphism(0,200)<439,0>[A`A/U;\pi_U]
\morphism(439,200)<382,0>[A/U`0;]
\morphism(-796,-200)<405,0>[0`\hat{\hat U};]
\morphism(-391,-200)<391,0>[\hat{\hat U}`\hat{\hat A};\hat{\hat\iota}_U]
\morphism(0,-200)<439,0>[\hat{\hat A}`\widehat{\widehat{A/U}};\hat{\hat\pi}_U]
\morphism(439,-200)<382,0>[\widehat{\widehat{A/U}}`0;]
\morphism(-391,200)<0,-400>[U`\hat{\hat U};\alpha_U]
\morphism(0,200)<0,-400>[A`\hat{\hat A};\alpha_A]
\morphism(439,200)<0,-400>[A/U`\widehat{\widehat{A/U}};\alpha_{A/U}]
\efig
\end{align}
The group $A/U$ is discrete, and thus $\alpha_{A/U}$ is a topological 
isomorphism by Theorem~\ref{thm:P:cd}. Therefore, (a) follows from the 
well-known Five Lemma for abelian groups.
In order to show (b), observe that $\alpha_U=\alpha_A\rest_U$ and 
$\alpha_U^{-1}=\alpha_A^{-1}\rest_U$, and that it suffices to check the 
continuity of a homomorphism on a neighborhood of $0$. Since $U$ is open, 
the statement follows.
\end{proof}

\begin{samepage}

\begin{ftheorem}[\pcite{2.6}{Chasco1}] \label{thm:compact:subgr}
Let $A \in \mathsf{Ab(Haus)}$ be such that $\alpha_A$ is injective, 
and let $K\leq A$ be a compact subgroup. Then:

\begin{list}{{\rm (\alph{enumi})}}
{\usecounter{enumi}\setlength{\labelwidth}{25pt}\setlength{\topsep}{-8pt}
\setlength{\itemsep}{-5pt} \setlength{\leftmargin}{25pt}}   

\item
$\alpha_{A/K}$ is injective;

\item
$\alpha_{A/K}$ is surjective if and only if $\alpha_A$ 
is so;

\item
$\alpha_{A/K}$ is an isomorphism of topological groups if and only if
$\alpha_A$ is so.

\end{list}
\end{ftheorem}

\end{samepage}

\begin{remark}
Although in Proposition~\ref{prop:compact:subgr} we assumed no separation 
axioms, here we do, because $\alpha_A$ being injective implies that 
$A$ is Hausdorff.
\end{remark}

\begin{proof} Since $\alpha_K$ is injective, by 
Proposition~\ref{prop:compact:subgr}, the exact sequence in 
(\ref{eq:compact:subgr}) gives rise to an exact sequence
\begin{align}
0\to\widehat{A/K}\to^{\hat\pi_K}
&\hat A\to<285>^{\hat\iota_K}\hat K\to<285> 0
\end{align}
where $\widehat{A/K}$ is an open subgroup of $\hat A$ and $\hat\iota_K$ is 
a quotient. Therefore, by Proposition~\ref{prop:open:subgp}, the lower row 
of the commutative diagram below is exact:
\begin{align}
\bfig
\morphism(-796,200)<405,0>[0`K;]
\morphism(-391,200)<391,0>[K`A;\iota_K]
\morphism(0,200)<439,0>[A`A/K;\pi_K]
\morphism(439,200)<382,0>[A/K`0;]
\morphism(-796,-200)<405,0>[0`\hat{\hat K};]
\morphism(-391,-200)<391,0>[\hat{\hat K}`\hat{\hat A};\hat{\hat\iota}_K]
\morphism(0,-200)<439,0>[\hat{\hat A}`\widehat{\widehat{A/K}};\hat{\hat\pi}_K]
\morphism(439,-200)<382,0>[\widehat{\widehat{A/K}}`0;]
\morphism(-391,200)<0,-400>[K`\hat{\hat K};\alpha_K]
\morphism(0,200)<0,-400>[A`\hat{\hat A};\alpha_A]
\morphism(439,200)<0,-400>[A/K`\widehat{\widehat{A/K}};\alpha_{A/K}]
\efig
\end{align}
Since $K$ is compact Hausdorff, $\alpha_K$ is an isomorphism of 
topological groups (cf. Proposition~\ref{thm:P:cd}), and thus (a) and (b) 
follow from the Five Lemma for abelian groups.

(c) If $\alpha_A$ is an isomorphism of topological groups, then 
$\alpha_{A/K}$ is a bijection by (a) and (b), and  the topologies of 
$A/K$ and $\widehat{\widehat{A/K}}$ coincide, because 
$\widehat{\widehat{A/K}}\cong \hat{\hat A}/\hat{\hat K}\cong A/K$.

Conversely, suppose that $\alpha_{A/K}$ is a topological isomorphism. Then
$\alpha_A$ is bijective by (a) and (b). We show that $\alpha_A$ is 
continuous. To that end, let $\mathcal{F}$ be a filter converging to 
$0$ in $A$. Then $\pi_K(\mathcal{F})\longrightarrow K$, and so
$\hat{\hat\pi}_K(\alpha_A(\mathcal{F}))=
\alpha_{A/K}(\pi_K(\mathcal{F}))\longrightarrow \hat{\hat K}$. By 
Lemma~\ref{lemma:comp:proper-new}(a),
the filter $\alpha_A(\mathcal{F})$ has a cluster point 
$x \in \hat{\hat K}\subseteq \hat{\hat A}$, because $\hat{\hat K}$ is 
compact. So, let $\mathcal{G} \supseteq \mathcal{F}$ be a filter such that
$\alpha_A(\mathcal{G})\longrightarrow x$. For $\chi\in\hat A$, one has
\begin{align}
x(\chi)=\lim \mathcal{\alpha_A(G)}(\chi) = \lim \chi(\mathcal{G}) =
\chi (\lim\mathcal{G})= 0
\end{align}
because $\chi$ is continuous and $\mathcal{G}\longrightarrow 0$, and 
thus $x=0$. Therefore, $0$ is the unique cluster point of 
$\alpha_A(\mathcal{F})$ for
every filter $\mathcal{F} \longrightarrow 0$. Hence, 
$\alpha_A(\mathcal{F})\longrightarrow 0$ for every $\mathcal{F} 
\longrightarrow 0$, and $\alpha_A$ is continuous.

Since the situation is completely symmetric (one could invert
$\alpha_K$, $\alpha_A$, and $\alpha_{A/K}$), the proof of the continuity 
of $\alpha_A^{-1}$ is similar.
\end{proof}

\subsection{Locally compact groups} We conclude this chapter with making a 
step toward the classical Pontryagin duality for locally compact Hausdorff 
abelian (LCA) groups. We denote by $\mathsf{LC}$ the category of locally 
compact Hausdorff spaces, and thus $\mathsf{Ab(LC)}$ is the category of 
LCA groups.

\begin{samepage}
\begin{proposition} \label{prop:LC:basic}
Let $A \in \mathsf{Ab(LC)}$. Then:

\begin{list}{{\rm (\alph{enumi})}}   
{\usecounter{enumi}\setlength{\labelwidth}{25pt}\setlength{\topsep}{-8pt}
\setlength{\itemsep}{-5pt} \setlength{\leftmargin}{25pt}}

\item
$\hat A \in \mathsf{Ab(LC)}$;

\item
$\alpha_A$ is continuous;

\item
$\alpha_A$ is injective.
\end{list}
\end{proposition}
\end{samepage}

\begin{proof}
(a) Clearly, $\hat A$ is Hausdorff. Let $U \in \mathcal{N}(A)$ such that  
$\bar{U}$ is compact. By  Proposition~\ref{prop:ev:cont}(a), 
$U^\vartriangleright$ is compact in $\hat A$, and by 
Proposition~\ref{prop:pol:basic}(a), $\bar U^\vartriangleright$ is a basic 
neighborhood of $0$ in $\hat A$. Since
$\bar U^\vartriangleright \subseteq U^\vartriangleright$, this completes 
the proof.

(b) Let $V \in \mathcal{N}(A)$ be such that $\bar V$ is compact.
By Theorem~\ref{thm:ev:kcont}, $\alpha_A$ is $k$-continuous, and so
$\alpha_A\rest_{\bar V}$ is continuous. In particular, $\alpha_A$ is 
continuous at $0$, and therefore it is continuous.

(c) Let $a \in A$ be a non-zero element.

First, suppose that $A$ is generated by $V \in \mathcal{N}(A)$ such 
that $\bar V$  is compact.  By 
replacing $V$ with $V \cup (-V) \cup \{a,-a\}$ if necessary, we may assume 
that $a \in V$ and $-V=V$ from the outset. By \cite[Lemma 2.42]{Rudin}, 
there is a  subgroup $H \leq A$ such that $H \cap V=\{0\}$, 
$H \cong \mathbb{Z}^n$ for some $n \in \mathbb{N}$, and $A/H$ is compact. 
The image $a+H$ of $a$ in $A/H$ is non-zero, because $a \in V$ and $V\cap 
H =\{0\}$. Thus, by Theorem~\ref{thm:P:cd}(b), $\alpha_{A/H}(a+H)\neq 0$, 
and so there is $\bar \chi \in \widehat{A/H}$ such that 
$\bar\chi(a+H)\neq 0$. Therefore, one has
$\alpha_A(a)(\chi)=\chi(a)=\bar\chi(a+H) \neq 0$, where
$\chi=\hat\pi_{H}(\bar \chi)  \in \hat A$. 
Hence, $\alpha_A(a)\neq 0$.

In the general case, pick $V \in \mathcal{N}(A)$ such that $-V=V$, 
$a \in V$, and $\bar V$ is compact, and put $U=\nolinebreak\langle V\rangle$, 
the subgroup generated by $V$. 
Then $U+V \subseteq U$, and thus $\bar U \subseteq \operatorname{Int} U$ 
(by Proposition~\ref{prop:sep:topgrp}(a)). Therefore, $U$ is an open 
subgroup of $A$. By what we proved so far, $\alpha_U$ is injective,  
Hence, by Theorem~\ref{thm:open:subgr}(a), $\alpha_A$ is injective. 
\end{proof}

\begin{remark}
The proof of (c) falls short of being self-contained because it is based 
on a Lemma from \cite{Rudin}. I hope to find a simple proof of the 
relevant parts of the Lemma. Suggestions???
\end{remark}


\begin{lemma} \label{lemma:LC:Gdelta}
Let $A \in \mathsf{Ab(LC)}$, and suppose that there is a countable family 
$\{U_n\} \subseteq \mathcal{N}(A)$ of neighborhoods of $0$ in $A$ such 
that $\bigcap U_n = \{0\}$. Then $A$ is metrizable.
\end{lemma}

\begin{proof}
By replacing each $U_n$ with $U_n^\prime \in \mathcal{N}(A)$ such that 
$\bar U_n^\prime \subseteq U_n$ and $\bar U_n^\prime$ is compact 
if necessary, we may assume that $\bigcap \bar U_n = \{0\}$, and 
$\bar U_1$ is compact (such $U_n^\prime$ exists, because $A$ is regular).
Furthermore, without loss of generality, we may assume that 
$U_{n+1} \subseteq U_n$. By Theorem~\ref{thm:Icount:met}, it suffices to 
show that $\{U_n\}$ is a base at $0$ for $A$. To that end, let $V \in 
\mathcal{N}(A)$. Then $\bigcap (\bar U_n\backslash \operatorname{Int} V) = 
\emptyset$, and $\{\bar U_n\backslash \operatorname{Int} V\}$ is a
decreasing family of closed subsets of the compact space $\bar U_1$. 
Therefore, there is $n_0$ such  that $\bar U_{n_0}\backslash 
\operatorname{Int} V =\emptyset$. Hence,
$\bar U_{n_0} \subseteq \operatorname{Int} V$, as desired.
\end{proof}

\begin{lemma} \label{lemma:LC:almost-met}
Let $A \in \mathsf{Ab(LC)}$, and let $V \in \mathcal{N}(A)$ be such that
$\bar V$ is compact. Then $V$ contains a compact subgroup $K$ of $A$ such 
that $A/K$ is metrizable.
\end{lemma}

\begin{proof} We construct a family $\{V_n\}$ of neighborhoods of $0$ in 
$A$. Put $V_1=(-V)\cap V$, and for $n \in \mathbb{N}$, 
using continuity of addition in $A$, pick  $V_{n+1} \in \mathcal{N}(A)$  
such that $V_{n+1}+V_{n+1} \subseteq V_n$ and $-V_{n+1}=V_{n+1}$.
Put $K = \bigcap\limits_{n=1}^\infty V_n$; $K$ is a subgroup of $A$, 
because $K - K \subseteq V_{n+1} - V_{n+1} \subseteq V_n$ implies
$K - K \subseteq K$. One has 
$K + V_{n+1}\subseteq V_{n+1} + V_{n+1} \subseteq V_n$, and thus
$\bar K \subseteq \operatorname{Int} V_n$ (cf. 
Proposition~\ref{prop:sep:topgrp}(a)). Therefore, 
$\bar K \subseteq K \subseteq \bar V$, and $K$ is compact.
The inclusion $K + V_{n+1}\subseteq V_n$ also yields
$K = \bigcap\limits_{n=1}^\infty (K+V_n)$, which means that
$\bigcap\limits_{n=1}^\infty \pi_K(V_n)  = K$ in $A/K$.
Hence, $\{\pi_K(V_n)\}$ satisfies the conditions of 
Lemma~\ref{lemma:LC:Gdelta}, and $A/K$ is metrizable, as desired.
\end{proof}

\begin{theorem}
The following statements are equivalent:

\begin{list}{{\rm (\roman{enumii})}}
{\usecounter{enumii}\setlength{\labelwidth}{25pt}\setlength{\topsep}{-8pt}
\setlength{\itemsep}{-4pt} \setlength{\leftmargin}{25pt}}

\item
For every second countable  $A \in \mathsf{Ab(LC)}$, $\alpha_A$ is an 
isomorphism of topological groups.

\item
For every $A \in \mathsf{Ab(LC)}$, $\alpha_A$ is an isomorphism of 
topological groups.

\end{list}
\end{theorem}

\begin{proof} (i) $\Rightarrow$ (ii): Let $A \in \mathsf{Ab(LC)}$. First,
suppose that $A$ is generated by $V \in \mathcal{N}(A)$ such that $\bar V$ 
is compact. By replacing $V$ with $-V \cup  V$ if necessary, we may assume 
that $-V = V$ from outset. Then $A=\bigcup\limits_{n=1}^\infty
\underbrace{\bar V + \cdots + \bar V}_\text{n times}$, so $A$ is 
$\sigma$-compact, and in particular, $A$ is Lindel\"of. Let $K \subseteq V$ 
be a compact subgroup such that $A/K$ is metrizable (cf. 
Lemma~\ref{lemma:LC:almost-met}). The quotient $A/K$ is second countable, 
because it is a continuous image of the Lindel\"of space $A$. (The 
properties Lindel\"of, separable, and second countable are 
equivalent for metrizable spaces.) Thus, $\alpha_{A/K}$ is an isomorphism 
of topological groups by our assumption. By 
Proposition~\ref{prop:LC:basic}(c), $\alpha_A$ is injective, and 
therefore, by Theorem~\ref{thm:compact:subgr}(c), $\alpha_A$ is an 
isomorphism of topological groups.

In the general case, pick $V \in \mathcal{N}(A)$ such that $\bar V$ is 
compact, and put $U=\langle V \rangle$, the subgroup generated by $V$. 
Then $U+V \subseteq U$, and thus $\bar U \subseteq \operatorname{Int} U$ 
(by Proposition~\ref{prop:sep:topgrp}(a)). Therefore, $U$ is an open 
subgroup of $A$. By what we proved so far, $\alpha_U$ is an isomorphism of 
topological groups, because $U$ is generated by a compact neighborhood 
of $0$. Hence, by Theorem~\ref{thm:open:subgr}(b), $\alpha_A$ is an 
isomorphism of topological groups.
\end{proof}

\include{GrpObj}

\appendix

\chapter*{Appendix}

\addcontentsline{toc}{chapter}{\protect {}Appendix}

\thispagestyle{empty}

\makeatletter
  \@addtoreset{subsection}{section}
  \renewcommand{\thesection}{\Alph{section}}
  \renewcommand{\thesubsection}{\thesection.\arabic{subsection}}
\makeatother

\section{Separation properties of topological groups}

\label{app:septop}

For $G \in \mathsf{Grp(Top)}$, one puts $\mathcal{N}(G)$ for the 
collection of neighborhoods of $e$ in $G$.

\begin{fproposition} \label{prop:sep:topgrp}
Let $G \in \mathsf{Grp(Top)}$. 

\begin{list}{{\rm (\alph{enumi})}}
{\usecounter{enumi}\setlength{\labelwidth}{25pt}\setlength{\topsep}{-8pt}
\setlength{\itemsep}{-5pt} \setlength{\leftmargin}{25pt}}

\item
If $VW \subseteq U$ for $U,V,W \in \mathcal{N}(G)$, then
$\bar V,\bar W \subseteq \operatorname{Int U}$.

\item
$G$ is regular.

\item
The following statements are equivalent:

\begin{list}{{\rm (\roman{enumii})}}
{\usecounter{enumii}\setlength{\labelwidth}{25pt}\setlength{\topsep}{-8pt}
\setlength{\itemsep}{-4pt} \setlength{\leftmargin}{25pt}}

\item
$G$ is $T_3$ (i.e., $T_1$ and regular);

\item
$G$ is Hausdorff;

\item
$G$ is $T_1$;

\item
$G$ is $T_0$.

\end{list}

\end{list}
\end{fproposition}

\begin{proof}
(a) Set $W_1=(\operatorname{Int} W)^{-1}$, and let $x \in \bar V$. Then
$W_1 \in \mathcal{N}(G)$ (because the group inversion is continuous), and 
$V \cap xW_1 \neq \emptyset$. Thus, there is 
$v \in W$ and $w_1 \in W_1$ such that $xw_1=v$, or $x=v w_1^{-1}$.
Therefore, $x \in V W_1^{-1} = V (\operatorname{Int} W)$. 
Since $V (\operatorname{Int} W)$ is open, it is contained in 
$\operatorname{Int} (VW) \subseteq \operatorname{Int} U$. Hence,
$x \in \operatorname{Int} U$, as desired.

(b) By continuity of the multiplication group $m\colon G \times G 
\rightarrow G$, for every $U \in \mathcal{N}(G)$ there is $V \in 
\mathcal{N}(G)$ such that $VV\subseteq U$, and thus  $\bar V 
\subseteq \operatorname{Int}U$ by (a).

(c) Implications (i) $\Rightarrow$ (ii)  $\Rightarrow$ (iii)
$\Rightarrow$ (iv) are obvious. In order to complete the proof, observe 
that by (b), $G$ is regular, and every $T_0$ regular topological space $X$ is 
$T_1$. (Indeed, if $x,y \in X$ are  distinct points, then one of them, say 
$x$, has a neighborhood $U_0$ such that $y \not\in U_0$. Then 
$F=X\backslash \operatorname{Int} U_0$ is a closed subset containing $y$ 
that does not contain $x$. Since $X$ is regular, there are disjoint open 
subsets $U$ and $V$ of $X$ such that $F \subseteq U $ and $x \in V$. This 
completes the proof.)
\end{proof}


\begin{fproposition} \label{prop:max:T2}
Let $G \in \mathsf{Grp(Top)}$, and set $N_G=\bigcap \mathcal{N}(G)$. Then:

\begin{list}{{\rm (\alph{enumi})}}
{\usecounter{enumi}\setlength{\labelwidth}{25pt}\setlength{\topsep}{-8pt} 
\setlength{\itemsep}{-5pt} \setlength{\leftmargin}{25pt}}

\item
$N_G$ is a compact normal subgroup of $G$;

\item
$G/N_G$ is Hausdorff;

\item
for every continuous homomorphism $\varphi\colon G \rightarrow H$ into a 
Hausdorff group $H$, $N_G \subseteq \ker \varphi$. In particular,  
$\varphi$ factors uniquely through $G \rightarrow G/N_G$.

\end{list}
\end{fproposition}

\begin{proof}
(a)  Let $U \in \mathcal{N}(G)$ and $g \in G$. By continuity 
of the group operations, there is $V \in \mathcal{N}(G)$ such that
$V V^{-1} \subseteq U$ and $g^{-1} V g \subseteq U$. Thus, 
$N_G N_G^{-1}\subseteq VV^{-1}\subseteq U$ and $g^{-1} N_G g \subseteq
g^{-1} V g \subseteq U$.
This is true for every $U \in \mathcal{N}(G)$, and therefore
$N_G N_G^{-1} \subseteq N_G$ and $g^{-1} N_G g \subseteq  N_G$.  Hence, 
$N_G$ is a normal subgroup. Since $N_G$ is contained in every neighborhood 
of $e$, its compactness is clear.

(b) For every $U \in \mathcal{N}(G)$, there is $V\in \mathcal{N}(G)$ such 
that $VV \subseteq U$, and thus $VN_G\subseteq V V \subseteq U$. 
Thus, $N_G \subseteq \bigcap \{VN_G\mid V \in \mathcal{N}(G)\} 
\subseteq \bigcap \mathcal{N}(G)=N_G$. Therefore, $G/N_G$ is $T_0$, and 
hence it is Hausdorff, by Proposition~\ref{prop:sep:topgrp}.

(c) If $H$ is Hausdorff, then  $\{e_H\}=\bigcap\mathcal{N}(H)$, and so
\begin{align}
\ker\varphi = \bigcap \{\varphi^{-1}(W) \mid W \in \mathcal{N}(H)\}
\supseteq \bigcap \mathcal{N}(G) = N_G,
\end{align}
as desired.
\end{proof}

The group $G/N_G$ is the {\em maximal Hausdorff quotient} of $G$. In 
categorical language, this means that $G/N_G$ is the reflection of 
$G \in \mathsf{Grp(Top)}$ in $\mathsf{Grp(Haus)}$.

\begin{flemma} \label{lemma:comp:proper-new}
Let $G \in \mathsf{Grp(Top)}$, $K\leq G$ be a compact subgroup, and put
$\pi_K\colon G \rightarrow G/K$ for the canonical projection.

\begin{list}{{\rm (\alph{enumi})}}
{\usecounter{enumi}\setlength{\labelwidth}{25pt}\setlength{\topsep}{-8pt} 
\setlength{\itemsep}{-2pt} \setlength{\leftmargin}{25pt}}

\item
If $\mathcal{F}$ is a filter in $G$ such that $K$ is a 
cluster point of $\pi_K(\mathcal{F})$ in $G/K$, then $\mathcal{F}$ has a 
cluster point in $K$.

\item
If $L \subseteq G/K$ is compact, then so is $\pi_K^{-1}(K)$.

\end{list}
\end{flemma}

\begin{proof}
(a) Assume that $\mathcal{F}$ has no cluster point in $K$.  Then each 
point $x \in K$ has a neighborhood $U_x$ such that $U_x$ does not mesh 
$\mathcal{F}$, that is, there is $F_x \in \mathcal{F}$ such that 
$U_x \cap F_x=\emptyset$. The collection $\{U_x\}_{x\in K}$ is an open 
cover of $K$, and thus it has a finite subcover $U_{x_1},\ldots,U_{x_l}$, 
because $K$ is compact. Set $U=U_{x_1}\cup\ldots\cup U_{x_l}$ and
$F=F_{x_1}\cap\ldots\cap F_{x_l}$. Since $\mathcal{F}$ is a filter,
$F \in \mathcal{F}$, and one has $U\cap F=\emptyset$ and $K \subseteq U$.
Using the group multiplication $m\colon G \times G \rightarrow G$,
the last inclusion can be expressed as $\{e\}\times K\subseteq m^{-1}(U)$. 
Since $m$ is continuous, $m^{-1}(U)$ is open, and thus,
by a tube-lemma type argument, there is $V \in \mathcal{N}(G)$ 
such that $V \times K \subseteq m^{-1}(U)$. In other words,
$VK \subseteq U$. Therefore, $VK \cap F =\emptyset$, and so
$VK \cap FK =\emptyset$. Hence, $\pi_K(V)\cap \pi_K(F)=\emptyset$.
Since $\pi_K(V)$ is a neighborhood of $K$ in $G/K$, this contradicts
the assumption that $K$ is a cluster point of $\pi_K(\mathcal{F})$.

(b) Let $\mathcal{F}$ be a filter meshing $\pi_K^{-1}(L)$ (i.e., $F 
\cap\pi_K^{-1}(L)\neq \emptyset$ for every $F \in\mathcal{F}$). Then 
$\pi_K(\mathcal{F})$ is a filter that meshes $L$, and thus has a cluster 
point $xK$, because $L$ is compact. Equivalently, $K$ is a cluster point 
of $\pi_K(x^{-1}\mathcal{F})$, and therefore $x^{-1}\mathcal{F}$ has a 
cluster point in $K$ (by (a)).  Hence, $\mathcal{F}$ has a cluster point in 
$xK$, which shows that $\pi_K^{-1}(L)$ is compact. 
\end{proof}

A combination of Proposition~\ref{prop:max:T2}(a) and 
Lemma~\ref{lemma:comp:proper-new} yields:

\begin{corollary} \label{cor:max:T2}
Compact subsets of $G/N_G$ are precisely the images of compact subsets of 
$G$. \qed
\end{corollary}

\section{Exponentiability and the compact-open topology}
\label{app:A}

Let $\mathfrak{T}$ be a full subcategory of $\mathsf{Top}$, the category 
of topological spaces and their continuous maps. Suppose that 
$\mathfrak{T}$ has finite  products; they may differ from the products in 
$\mathsf{Top}$, but each product in $\mathfrak{T}$ must have the same 
underlying set as in $\mathsf{Top}$. 
Each function $f\colon X \times Y \rightarrow Z$  
gives rise to a map $\bar f \colon Y \rightarrow Z^X$ 
defined by  $\bar f(y) = f_y$, where $f_y \in Z^X$ is given by 
$f_y(x) = f(x,y)$.  One would like to equip the function set 
$\mathfrak{T}(X,Z) \subseteq Z^X$
with a suitable topology such that $\bar f$ is continuous whenever $f$ is 
so. If we add the natural requirement of functoriality of the topology, 
we arrive at the concept of  exponentiability. A space $X \in 
\mathfrak{T}$ is {\em exponentiable} (in $\mathfrak{T}$) if the functor
$X \times - : \mathfrak{T} \longrightarrow \mathfrak{T}$
has a right adjoint.  Since we assume $\mathfrak{T}$ to have finite 
products, it has a terminal object which must then be the one-point space 
$\{*\}$, because $\mathfrak{T}$ is a full subcategory. Thus, one expects 
$\mathfrak{T}(X,Z)$  to be the underlying set of the hom-object. In 
other words, $X$ is exponentiable if there is a  way of topologizing  the 
set $\mathfrak{T}(X,Z)$  for every $Z \in \mathfrak{T}$ such that the 
resulting space is in $\mathfrak{T}$, and the map
\begin{align}
\mathfrak{T}(X \times Y,Z) & \longrightarrow 
\mathfrak{T}(Y,\mathfrak{T}(X,Z)) 
\label{eq:kspaces:exp1}\\
f &  \longmapsto \bar f \nonumber
\end{align}
is a bijection that is natural in $Y$ and $Z$.

The question of describing exponentiable Hausdorff spaces was raised by 
Hurewicz in a personal communication with Fox, who was the first to give a 
partial answer to the question (cf. \cite{Fox}). He proved that for every 
regular locally compact space $X$ and every $Y \in \mathsf{Haus}$ the set
$\mathsf{Haus}(X,Y)$ can be topologized in a way that $f$ is continuous if 
and  only if $\bar f$ is so. Fox also showed that for every 
separable metric space $X$, it is possible to topologize
$\mathsf{Haus}(X,\mathbb{R})$ in a way that (\ref{eq:kspaces:exp1}) holds 
for every $Y$ if and only if $X$ is locally compact. Although Fox did not 
prove that the condition of local compactness is sufficient, and he showed 
necessity only for separable metric spaces $X$, he was actually very close 
to a complete solution.

For Hausdorff spaces $X$ and $Y$, we denote by $\mathscr{C}(X,Y)$ the 
space $\mathsf{Haus}(X,Y)$ equipped with the {\em compact-open topology}: 
Its subbase is the family
$\{[K,U] \mid K \subseteq X \text{ compact}, U \subseteq Y\text{ open}\}$,
where $[K,U] = \{ f \in \mathsf{Haus}(X,Y) \mid f(K) \subseteq U\}$.
The compact-open topology was also invented by Fox, and 
following his work, Arens studied the separation properties 
of $\mathscr{C}(X,Y)$ (cf. \cite{Arens}). Arens  was not 
far from proving that local compactness is necessary for exponentiability 
in $\mathsf{Haus}$ (cf. \cite[\nolinebreak Theorem~3]{Arens}). Jackson 
proved that if $X$ is locally 
compact, then  $\mathscr{C}(X,-)$ is the right adjoint of $X \times -$, 
moreover, the bijection
\begin{equation} \label{eq:kspaces:Jackson}
\mathscr{C}(X \times Y,Z) \longleftrightarrow 
\mathscr{C}(Y,\mathscr{C}(X,Z))
\end{equation}
is actually a homeomorphism for every Hausdorff space $Y$ and $Z$ (cf.  
\cite{Jackson}). It appears that Fox already recognized that 
(\ref{eq:kspaces:Jackson}) is a bijection, and Jackson's main 
achievement is proving that it is a homeomorphism.

\begin{remark}
Some authors, including Isbell , erroneously credit Brown for proving 
(\ref{eq:kspaces:Jackson}) (cf. \cite{Isbell1}). The fact is that Brown 
himself refers both to Fox and Jackson in his paper in question, and he 
lays no claim to the ``classical" results on special cases of the 
exponential law (as a personal communication with him reveals). 
At the same time, Brown was the first to show that the 
category of Hausdorff $k$-spaces is cartesian closed, but 
unfortunately this seems to be somewhat forgotten (cf.~\cite[3.3]{Brown}).
\end{remark}

Whitehead proved that if $X$ is a locally compact Hausdorff space, then 
$1_X \times g$ is a quotient map in $\mathsf{Haus}$ for every quotient map 
$g$ in $\mathsf{Haus}$ (cf. \cite{Whitehead}). Michael proved that the 
converse is also true: For a Hausdorff space $X$, the map $1_X \times g$ is 
a quotient map in $\mathsf{Haus}$ for every quotient map $g$ in 
$\mathsf{Haus}$ if and only if $X$ is locally compact 
(cf.~\cite[2.1]{Mich0}). In other words, the functor $X \times -$ 
preserves quotients if and only if $X$ is locally compact. If $X$ is 
exponentiable, then $X \times -$ must preserve every colimit in 
$\mathsf{Haus}$ (\cite[V.5.1]{MacLane}), and in particular it has to 
preserve coequalizers in $\mathsf{Haus}$ (which is equivalent to 
preservation of quotients; for details see~\cite[1.1]{GLPHD}). Thus, by 
Michael's result the local compactness of $X$ follows. Therefore, 
exponentiable spaces in $\mathsf{Haus}$ can be characterized as follows.

\begin{ftheorem} \label{thm:LC:exp}
A space $X \in \mathsf{Haus}$ is exponentiable in $\mathsf{Haus}$ if and 
only if $X$ is locally compact.
\end{ftheorem}

A proof of sufficiency of local compactness in Theorem~\ref{thm:LC:exp}
is available in standard topology textbooks (cf.~\cite[3.2]{Engel6}).

\begin{fcorollary} \label{cor:CO:kcont}
Let $X,Y \in \mathsf{Haus}$. The evaluation map
\begin{align}
e: X & \longrightarrow \mathscr{C}(\mathscr{C}(X,Y),Y) \\
x & \longmapsto \hat x \quad  [\hat x(f)=f(x)]
\end{align}
is $k$-continuous, that is, $e\rest_K$
is continuous for every compact subset $K \subseteq X$.
\end{fcorollary}

\begin{proof}
Let  $K\subseteq X$ be a compact subset. 
The map
\begin{align}
E\colon K \times \mathscr{C}(X,Y) & \longrightarrow Y \\
(x,f)&\longmapsto f(x)
\end{align}
is continuous, because $E(K\times [K,W]) \subseteq W$ for every open
subset $W$ of $Y$. Let $[\Phi,U]$ be a subbasic open subset of
$\mathscr{C}(\mathscr{C}(X,Y),Y)$, where $\Phi \subseteq \mathscr{C}(X,Y)$ is 
compact and $U \subseteq Y$ is open. By Theorem~\ref{thm:LC:exp}, $\Phi$
is exponentiable (because it is compact), and thus 
$E\rest_{K\times\Phi}\colon K \times \Phi
\rightarrow Y$  gives rise to a continuous map 
$e_{K,\Phi}\colon K \rightarrow \mathscr{C}(\Phi,Y)$ that is 
defined by $(e_{K,\Phi}(k))(f)=f(k)$. Therefore, 
$(e\rest_K)^{-1}([\Phi,U])=e_{K,\Phi}^{-1}([\Phi,U])$ is open in $K$, as 
desired.
\end{proof}

Exponentiable spaces in $\mathsf{Top}$ were characterized by Day and Kelly 
as so-called {\em core-compact} spaces (cf.~\cite{DayKel}), but this 
result is beyond the scope of this summary. For a detailed review of 
exponentiability, we refer the reader to \cite[Chapter 1]{GLPHD}.

\section{Hausdorff {\itshape k}-spaces}

\label{app:kspace}

A map $f\colon X \rightarrow Y$ between Hausdorff spaces is said to be 
{\em $k$-continuous} if the restriction $f\rest_K$ is continuous for every 
compact subset $K$ of $X$. Theorem~\ref{thm:ev:kcont} suggests that it 
might be helpful to consider spaces whose $k$-continuous maps are 
continuous.

\begin{samepage}

\begin{fproposition}[\pcite{3.3.18-21}{Engel6}] \label{prop:k:def}
Let $X \in \mathsf{Haus}$. The following statements are equivalent:

\begin{list}{{\rm (\roman{enumi})}}
{\usecounter{enumi}\setlength{\labelwidth}{25pt}\setlength{\topsep}{-6pt}
\setlength{\itemsep}{-5pt} \setlength{\leftmargin}{25pt}}

\item
if $O \cap K$ is open in $K$ for every $K \in \mathcal{K}(X)$, then $O$ is
open in $X$;

\item
if $F \cap K$ is closed for every $K \in \mathcal{K}(X)$, then $F$ is
closed;

\item
every $k$-continuous map $f\colon X \rightarrow Y$ into a Hausdorff
space $Y$ is continuous.

\end{list}
\end{fproposition}

\end{samepage}

\begin{proof}
The equivalence of (i) and (ii) is obvious, and clearly each of them
implies (iii).

(iii) $\Rightarrow$ (i): Let $\mathsf k X$ be the underlying set of $X$
equipped with the topology
\begin{align}
\{O \subseteq X\mid O\cap K\mbox{ is open in $K$ for every }
K\in\mathcal{K}(X)\}.
\end{align}
We call $\mathsf{k}X$ the {\em $k$-ification} of $X$.
Its topology is finer than the topology of $X$, and so
$\mathsf{k} X$ is Hausdorff. Furthermore, the topologies of $X$ and
$\mathsf{k} X$ coincide on compact subsets, and thus the identity map
$X \rightarrow \mathsf{k} X$ is $k$-continuous. Therefore, by (iii),
it is continuous, and hence $X \cong \mathsf{k} X$, as desired.
\end{proof}

A Hausdorff space $X$ is a {\em $k$-space} if it satisfies the equivalent 
conditions of Proposition~\ref{prop:k:def}. The category of Hausdorff 
$k$-spaces and their continuous maps is denoted by $\mathsf{kHaus}$. It is 
possible to define $k$-spaces without assuming any separation axioms 
(cf.~\cite{Vogt} and \cite[1.1]{GLPHD}), but it is unnecessary for our 
purposes.

\begin{ftheorem} \mbox{ }

\begin{list}{{\rm (\alph{enumi})}}
{\usecounter{enumi}\setlength{\labelwidth}{25pt}\setlength{\topsep}{-6pt}
\setlength{\itemsep}{-5pt} \setlength{\leftmargin}{25pt}}

\item
The category $\mathsf{kHaus}$ is a coreflective subcategory of 
$\mathsf{Haus}$
with coreflector $\mathsf{k}$, the $k$-ification.

\item
The product of $X, Y \in \mathsf{kHaus}$ in $\mathsf{kHaus}$ is the 
$k$-ification $\mathsf{k}(X\times Y)$ of their product in $\mathsf{Haus}$.

\item {\rm (\cite[3.3]{Brown})}
The category $\mathsf{kHaus}$ is cartesian closed, and the internal hom
is given by the $k$-ification $\mathsf{k}\mathscr{C}(X,Y)$ of the 
compact-open topology.

\end{list}
\end{ftheorem}

\begin{proof}
(a) First, note that by Proposition~\ref{prop:k:def}, $\mathsf{k}X$ is a 
Hausdorff $k$-space for every $X \in \mathsf{Haus}$.
Let $g\colon X \rightarrow Y$ be a continuous map between Hausdorff 
spaces $X$ and $Y$, and let $O \subseteq \mathsf{k}Y$  be open. Since $g$ 
is continuous, $g(K)\in\mathcal{K}(Y)$ for every $K\in\mathcal{K}(X)$.
Thus, $g(K) \cap O$ is open in $g(K)$, and hence 
$K\cap g^{-1}(O) = g^{-1}(g(K) \cap O)$ is open in $K$. Therefore,
$\mathsf{k}g \colon \mathsf{k} X \rightarrow \mathsf{k} Y$ is continuous.
This shows that $\mathsf{k}\colon \mathsf{Haus}\longrightarrow 
\mathsf{kHaus}$ is a functor (it is clear that it preserves 
$\operatorname{id}_X$ and composition). In order to show that $\mathsf{k}$ 
is a coreflector, observe that if $X \in \mathsf{kHaus}$ (i.e., $X 
=\mathsf{k}X$), then $\mathsf{k} f\colon X \rightarrow \mathsf{k} Y$ is 
the unique continuous map such that $f=r_Y\circ \mathsf{k} f$, where
$r_Y\colon\mathsf{k}Y \rightarrow Y$  is the identity.

(b) follows from (a) and \cite[VI.3, V.5.1]{MacLane}.

It is not hard to see that (c) follows from Corollary~\ref{cor:CO:kcont}, 
but for a detailed proof,  we refer the reader to \cite[VII.8.3]{MacLane}.
\end{proof}

\bibliography{notes,notes2,notes3}

\end{document}